\documentclass[a4paper,11pt]{article}
\usepackage{hyperref}
\usepackage{geometry}
\usepackage{comment}
\usepackage{authblk}

\usepackage[british]{babel}
\usepackage{wrapfig}
\usepackage[comma, square, numbers]{natbib}
\usepackage{yfonts}
\usepackage[utf8]{inputenc}
\usepackage{amssymb}
\usepackage{amsthm}
\usepackage{graphics}
\usepackage{soul}
\usepackage{amsmath}
\usepackage{amstext}
\usepackage{engrec}
\usepackage{floatflt}
\usepackage{rotating}
\usepackage[safe,extra]{tipa}
\usepackage[font={footnotesize}]{caption}
\usepackage{framed}
\usepackage{listings}
\usepackage{subcaption}
\usepackage{enumerate}
\usepackage{tikz}

\usepackage{mathtools}

\usepackage[titletoc,title]{appendix}

\newcommand{\vertiii}[1]{{\left\vert\kern-0.25ex\left\vert\kern-0.25ex\left\vert #1 
    \right\vert\kern-0.25ex\right\vert\kern-0.25ex\right\vert}}

\newcommand{\mc}{\mathcal}
\newcommand{\mb}{\mathbb}

\newcommand{\R}{\mb R}

\newcommand{\N}{\mb N}
\newcommand{\Z}{\mb Z}

\newcommand{\T}{\mb T}

\newcommand{\eea}{\end{align}}

\renewcommand{\epsilon}{\varepsilon}
\renewcommand{\bar}{\overline}

\newcommand{\bo}{\boldsymbol}
\renewcommand{\phi}{\varphi}

\renewcommand\upsilon{\theta}

\usepackage[normalem]{ulem}

\newtheorem{theorem}{Theorem}[section]

\newtheorem{corollary}{Corollary}[section]
\newtheorem{lemma}{Lemma}[section]

\newtheorem{proposition}{Proposition}[section]
\theoremstyle{definition}
\newtheorem{definition}{Definition}[section]

\theoremstyle{remark}
\newtheorem{remark}{Remark}[section]

\newtheoremstyle{algorithm}
{4pt}
{4pt}
{}
{}
{}
{:}
{\newline}
{}

\usepackage{xcolor,placeins}

\newtheorem{algorithm}{Algorithm}

\newcommand{\balgorithm}{\begin{algorithm}\begin{framed}\ }
\newcommand{\ealgorithm}{\end{framed}\end{algorithm}}

\newcommand{\pippo}{\begin{code} \  \begin{lstlisting}[frame=single]}
\newcommand{\pappo}{\end{lstlisting} \end{code}}

\newcommand{\bd}{\begin{definition}}
\newcommand{\ed}{\end{definition}}

\newcommand{\bt}{\begin{theorem}}
\newcommand{\et}{\end{theorem}}
\newcommand{\bp}{\begin{proposition}}
\newcommand{\ep}{\end{proposition}}
\newcommand{\bc}{\begin{corollary}}
\newcommand{\ec}{\end{corollary}} 
\newcommand{\bl}{\begin{lemma}}
\newcommand{\el}{\end{lemma}}
\newcommand{\br}{\begin{remark}}
\newcommand{\er}{\end{remark}}

\DeclareMathOperator{\End}{End}

\DeclareMathOperator{\Id}{Id}

\theoremstyle{definition}

\theoremstyle{remark}

\theoremstyle{example}

\title{Mean-Field Coupled Systems and Self-Consistent Transfer Operators: A Review}
\date{}
\author[1]{Matteo Tanzi\footnote{mtanzi@lpsm.paris}}
\affil[1]{Laboratoire Statistique Probabilité Modélisation, CNRS-Sorbonne Universit\'e-Universit\'e Paris Cit\'e.}
\begin{document}
\maketitle

\begin{abstract}
In this review we survey the literature on mean-field coupled maps. We start with the early works  from the physics literature, arriving to some recent results from ergodic theory  studying the thermodynamic limit of globally coupled maps and the associated self-consistent transfer operators. We also give few pointers to related research fields dealing with mean-field coupled systems in continuous time, and applications. 
\bigskip

\bigskip

\noindent
{\bf Statements and Declarations:} The author has no conflict of interests to declare.
\bigskip

\noindent
{\bf Acknowledgments:}  The author was supported by the Marie Sk{\l}odowska-Curie actions project: ``Ergodic Theory of Complex Systems" proj. no. 843880. 
\end{abstract}

\section{Introduction}
Understanding the dynamics of complex systems is  the forefront of  research in many areas of science. Examples of complex systems most impactful in our everyday lives are networks of neurons, gene regulatory networks, artificial neural networks, spread of epidemics, and opinion models. The behavior of these systems is the result of the intricate interactions of their microscopic components. 
 
Oftentimes, complex systems are modeled as dynamical systems interacting on a graph/network: each dynamical system represents a fundamental component of the complex system (e.g. a gene, a neuron, an individual) and occupies a node on the graph whose edges prescribe the interactions between components. The literature on dynamical systems coupled on networks is vast. The tutorial \cite{porter2016dynamical} presents a broad survey.

In this review we focus on mean-field coupled systems in discrete time (\emph{globally coupled maps}) and their thermodynamic limits (\emph{self-consistent transfer operators}).  Mean-field models are characterized by  a very large number of components coupled by weak pairwise interactions whose strength scales as the inverse of the number of coupled units.  Numerical simulations show that these systems exhibit a great variety of behaviors, some of which are reminiscent of complex systems.

Most rigorous arguments deal with  thermodynamic limits of  coupled systems where the number of components tends to infinity. In this limit, the global state is given by a probability measure describing the distribution of the infinitely many components in phase-space, and its evolution is given by a nonlinear evolution law prescribed by a self-consistent transfer operator.  

Below we  report the main observations and results available on globally coupled maps.  The literature on the topic is vast, and a complete review of all the contributions to the subject seems hopeless. Rather than to be complete, the objective of this paper is to: give some history and context to the study of mean-field coupled maps, survey the  advances in the study of self-consistent transfer operators made in the last decade, and provide pointers to  research topics  having an affinity with mean-field coupled maps. 

\smallskip
\emph{Organization of the review.} In Section \ref{Sec:GloballyCoupledMaps} we review the literature on globally coupled maps from the origins in the physics literature (sections \ref{Sec:GlobCouplVarBehav}-\ref{Sec-LYap}) and the ergodic theoretical approaches to study these systems in (Section \ref{Sec:ErgThCoupMaps}). We also review other coupled systems in discrete time, in particular maps coupled on lattices and heterogeneous networks (Section \ref{Sec:Others}). In Section \ref{Sec:STOs} we focus on the thermodynamic limits  and  the study of self-consistent transfer operators. We review various rigorous approaches to study existence and stability of fixed states (Section \ref{Sec:StabFixedPoints}) and their  linear response (Section \ref{Sec:LinResp}). We then look at situations where the self-consistent operators have more complicated attractors, and the available studies are mostly numerical (sections \ref{Sec:TwardsmreComp}-\ref{Sec:NumericalStudies}). We conclude the section reviewing a recent development on propagation of chaos for globally coupled maps. Finally in Section \ref{Sec:Othertopics} we give some pointers to works dealing with mean-field models in continuous time and other related topics. Among others, we give very quick (and superficial) overviews of: interacting particle systems, systems  of coupled oscillators, mean-field models on adaptive and higher-order networks, highlighting the connections with globally coupled maps. 

\section{Globally  Coupled Maps}\label{Sec:GloballyCoupledMaps}
At the end of the '80s beginning of the '90s, \emph{globally coupled maps} (GCMs) arose as high-dimensional models of complex systems having simple equations, but whose dynamics exhibited a great variety of behaviors.  Loosely speaking, the equations describing the evolution of $N$ identical coupled \emph{maps}, also called \emph{units} or \emph{sites}, have the form 
\begin{equation}\label{Eq:GloballyCoupled}
x_i(t+1)=f(x_i(t))+\frac{1}{N}\sum_{j=1}^Nh(x_i(t),x_j(t)),\quad\quad i=1,...,N
\end{equation}
where $x_i(t)$ characterizes the state at each site and belongs to $M$, a set with some additive operation  (often an interval or $\T=\R\backslash \Z$). The map $f$ is called the \emph{local} or \emph{uncoupled} dynamics. Each term $\frac1N h(x_i,x_j)$ gives the pairwise additive interaction that the $i$-th  unit receives from the $j$-th one. Different formulations can be found in the literature some of which will be discussed through this review.

This model stemmed from a similar setup where maps are coupled on a lattice: Given $d\ge 1$ and $\Lambda\subset \Z^d$ finite or infinite,
\begin{equation}\label{Eq:CML}
x_i(t)=f(x_i(t))+\epsilon\sum_{j\in\Lambda_i}h(x_i(t),x_j(t))\quad\quad i\in\Lambda
\end{equation}
where $\Lambda_i\subset\Lambda$ prescribes a set of neighbors of $i$. The above, is a continuous variable version of a spin system, and was introduced as a model to study chaos and pattern formation in spatially extended systems after coarse graining. Systems as in \eqref{Eq:CML} are called \emph{coupled map lattices} (for a brief review see Section \ref{SeubSec:CML}).   

Globally coupled maps are coupled map lattices where the set of neighbors $\Lambda_i$ is the whole collection of units, and where $\epsilon$ is required to scale as $N^{-1}$. One reason behind this normalization invokes energy considerations: If a unit has to spend some ``energy'' to influence another unit and has only a finite amount of ``energy'' to spend, this will have to be distributed among all the interactions made. If interactions are identical, it will have to be distributed equally. This assumption is chiefly made for  the equations to be well defined for  $N\rightarrow\infty$. If the interaction strength scales as $N^{-1}$, in the limit $N\rightarrow +\infty$ one hopes that the mean-field interaction term in \eqref{Eq:GloballyCoupled} converges to a number depending only on the global distributions of the $\{x_i(t)\}_{i=1}^N$, and not on their particular state, making the equation that defines the dynamic of $x_i(t)$ virtually only dependent on $x_i(t)$ and identical across units.


\subsection{Synchronization, phase ordering, and turbulence}\label{Sec:GlobCouplVarBehav}
It  does not come as a surprise that globally coupled maps exhibit a great range of behaviors that can vary changing the parameters that define the local dynamics and/or the coupling among units. The asymptotic behavior of the orbits can also drastically change depending on the initial condition, which suggests that these systems possess a large number of different attractors. 

\medskip
 One of the first instances where the above observations were reported is \cite{kaneko1989chaotic}\footnote{Similar observations in coupled oscillators had appeared in \cite{wiesenfeld1989attractor}.}. Here the globally coupled system considered is 
\begin{equation}\label{Eq:CoupledLogisticMaps}
x_i(t+1)=(1-\epsilon)f(x_i(t))+\frac{\epsilon}{N}\sum_{j=1}^Nf(x_j(t)),\quad\quad i=1,...,N
\end{equation}
with  $x_i(t)\in[-1,1]$ and $\epsilon\ge 0$, and corresponds to the application of  an uncoupled map $f:[-1,1]\rightarrow[-1,1]$, and of
\[
x_i(t+1)=(1-\epsilon)x_i(t)+\frac{\epsilon}{N}\sum_{j=1}^Nx_j(t) \quad\quad i=1,...,N
\]
 a \emph{diffusive} mean-field interaction where the state of each unit tends to get closer to the average of the states of all the units. Notice that for $\epsilon=0$ the maps are uncoupled and evolve according to $f$, while for $\epsilon=1$, after one time step the system synchronizes instantaneously and each coordinate takes the same value equal to the average state. In between these values a great variety of behaviors can arise.

In  \cite{kaneko1989chaotic}, $f$ belongs to the logistic family 
\begin{equation}\label{Eq:LogisiticFam}
f(x):=f_{a}(x)=1-ax^2.
\end{equation}
 which is known to exhibit intricate bifurcation patterns alternating periodic and chaotic attractors.
With this choice, varying $\epsilon\ge 0$ and $0\le a\le 2$ and simulating the dynamics for different initial conditions, various behaviors are observed \cite{kaneko1990clustering} which have been organized in terms of their synchronization patterns. Two units  $i$ and $j$ are \emph{synchronized} if $x_i(t)=x_j(t)$\footnote{Different notions of synchronization exist. We do not address these differences here, but we point the reader to \cite{pikovsky2002synchronization}. On synchronization, see also \cite{boccaletti2002synchronization} and \cite{arenas2008synchronization}.}. A cluster, is a maximal subset  of the units, $\{i_1,...,i_k\}\subset\{1,...,N\}$, whose units are synchronized, i.e.  $x_{i_1}(t)=...=x_{i_k}(t)$.

One   says that an orbit exhibits
\begin{enumerate}
\item \emph{Synchrony}: if there is only one cluster, i.e. all the units are synchronized;
\item \emph{Ordered phase}: if there is a ``small" number of clusters each of which contains a ``large" fraction of the units;
\item \emph{Partially ordered phase}: there is a large number of clusters with few of them containing a large fraction of the units, and many of them containing few units;
\item \emph{Turbulent/chaotic phase}: there is no discernible organization into clusters.
\end{enumerate}
For an $(a,\epsilon)$-bifurcation diagram showing the emergence of these different behaviors see e.g. Figure 1 in \cite{kaneko1989chaotic}. As expected, small values of $\epsilon$ favor the turbulent phase, while larger values of $\epsilon$ favor the emergence of coherent structures as in (1) or (2). One can classify states from (1)-(3) with respect to the number and size  of clusters by associating to the sate an $m$-tuple $(k_1,...,k_m)$ with $k_1+...+k_m=N$ that denote the presence of $m$ clusters with cluster $i$ containing $k_i$ units. 

\smallskip
For fixed parameters, one observes a wide  range  of dynamics  already  among states having the same number of clusters, only differing for the clusters' size. For example, if $m=2$ one can cook up the parameters so that  if $k_1\approx k_2\approx N/2$, the motion is  periodic with  the units in each cluster switching between two different values $X_{1*},\,X_{2*}$ at each time-step in antiphase, i.e. when cluster 1 is close to $X_{1*}$, cluster 2 is close to $X_{2*}$ and vice versa. However, moving units from cluster 2 to cluster 1 and thus increasing the value of $k_1$ at the expense of  $k_{2}$, one assists to a period doubling cascade with the states of the first cluster jumping periodically between $2^n$ values (see e.g. Figure 13a in \cite{kaneko1990clustering}). Further increase of $k_1$ destroys the two clusters and starts a transition to the turbulent regime. The above picture is explained by the following observations: Letting 
\[
X_1(t):=x_1(t)=...=x_{k_1}(t)\quad\mbox{ and }\quad X_2(t):=x_{k_1+1}(t)=...=x_N(t)
\] 
be the states of the units in each cluster, and substituting in \eqref{Eq:CoupledLogisticMaps},  the evolution equations for $X_1(t)$ and $X_2(t)$  are 
\begin{align*}
X_1(t+1)&=(1-\epsilon_2)f(X_1(t))+\epsilon_2 f(X_2(t))\\
X_2(t+1)&=(1-\epsilon_1)f(X_2(t))+\epsilon_1 f(X_1(t))
\end{align*}
with $\epsilon_1:=\epsilon\frac{k_1}{N}$ and $\epsilon_2:=\epsilon\frac{k_2}{N}$. Thus, changing $k_1$ produces a change in $\epsilon_1,\,\epsilon_2$ and a bifurcation in the 2D system above that can explain the observed  behavior. Let us stress that these bifurcations are for fixed values of $(a,\epsilon)$ varying $k_1$, i.e. they are observed in the \emph{same} system only changing  the initial condition.

For a general  account on the bifurcations with respect to the $(a,\epsilon)$ parameters see \cite{balmforth1999synchronized} and in particular Figure 1 therein which captures the variety of attractors numerically observed in the system.  

\smallskip
 All of the  above suggest that these systems have a large number of attractors. A similar situation was observed in \cite{wiesenfeld1989attractor} in a system of coupled continuous time oscillators where this coexistence  phenomenon has been termed \emph{attractor crowding}. An important feature of attractor crowding is that  attractors increase  factorially in number with the system size -- \cite{wiesenfeld1989attractor} estimates $(N-1)!$ -- and get closer in phase space so that a small perturbation of an orbit can drive the system from one attractor to the other giving high-versatility expected to have implications on the system's function (see also Sect. \ref{Sec:ApplicationofMeanField} below on applications).  
\medskip

Parameters $(a,\epsilon)$ can be chosen to give rise to the following interesting phenomenon called \emph{posi-nega switching} \cite{kaneko1989chaotic,kaneko1990clustering}: Start from an initial condition in a two-clusters state $(k_1,k_2)$ with $k_1\approx k_2$, where the dynamic of each cluster follows a periodic orbit of period 2 in antiphase. Now, assume that perturbing the state of a single map  in the second cluster of a fixed quantity $\delta$ brings it to a region where the map eventually joins the first cluster. This leads to an increase of $k_1$  and  a period doubling cascade (see Figure 2 in \cite{kaneko1989chaotic} and Figure 3 (c) in \cite{shimada2000resolution}) that roughly corresponds to the period doubling cascade in the logistic map. Further increase of $k_1$ makes each cluster split with the units separating  and  undergoing chaotic motion. At this point, there is no more division into clusters. However, if we keep adding the fixed perturbation to the maps that used to belong to the second cluster, at some point the two clusters reform and recover their periodic switching dynamic.  What is most surprising is that clusters reform so that maps that were in cluster 1 before the chaotic phase join the same cluster, and similarly for maps originally in cluster 2 suggesting that the system keeps some ``memory" of the cluster subdivision.

\medskip
In \cite{kaneko1991globally},  a similar picture to the one in \eqref{Eq:CoupledLogisticMaps}  has been presented with coupled maps having $x_i(t)\in[0,1]$ governed by the equations
\[
x_i(t+1)=x_i(t)+\frac{K}{N}\sum_{j=1}^N\sin[2\pi x_i(t)-2\pi x_j(t)],\quad\quad i=1,...,N.
\]
In contrast with the previous setups where the emergence of ordered, partially ordered phases, and the richness of periodic orbits stemmed from the interplay between the richness of the local dynamics and the diffusive coupling, here it is only due to the coupling.

\medskip
We also mention \cite{cosenza1998synchronization} which studies  logarithmic maps, $f_a(x)= a+\ln|x|$, coupled as in \eqref{Eq:CoupledLogisticMaps}. Here synchronization and collective behaviors are observed, but there is no subdivision of the units into clusters (there is only one cluster), and the state $x(t)=x_1(t)=x_2(t)=...$ of the synchronized units  either evolves around a periodic orbit, or undergoes  chaotic motion. For certain values of $(a,\epsilon)$, a turbulent phase is  detected. See Figure 4 in \cite{cosenza1998synchronization} for a bifurcation diagram.

\medskip
Some works focused on the case of uniformly hyperbolic local dynamics given by uniformly expanding maps or tent maps in the parameter regimes having a unique  mixing absolutely continuous invariant measure. In \cite{just1995globally}, the author studies  globally coupled system with $x_i(t)\in[0,2\pi]$ governed by the equations
\[
x_i(t+1)=f(x_i(t))+\frac{\epsilon}{N}\sum_{j=1}^N\sin[f(x_j(t))-f(x_i(t))] \mod 2\pi,\quad\quad i=1,...,N
\] 
where $f$ depends on a parameter $a\in\R$, and is given by
\[
f(x)=2x+a\sin(x).
\]
It is immediate to see that for $a=\epsilon=0$, the system is a product of $N$ uncoupled doubling maps, therefore it has a unique absolutely continuous invariant mixing probability measure on $[0,2\pi]^N$. Nowadays, it is well known that this picture is stable under small perturbations, and therefore persists for small $a$ and $\epsilon$. At the time, this was claimed in \cite{just1995globally} using Markov partitions. In the same paper, it is shown that larger values of $\epsilon$ lead to synchronization, i.e. the synchronization manifold 
\[
\left\{(x_1,...,x_N)\in[0,2\pi]^N:\quad x_1=...=x_N\right\}
\]  
is a stable invariant set. The above picture shows that the dynamics undergoes bifurcations shifting from a regime where the behavior is dictated by the hyperbolic local dynamics, to a situation where the behavior is dictated by the coupling. The presence of a unique absolutely continuous measure in the small coupling regime was rigorously proved in  \cite{keller1997mixing} for systems of coupled tent maps. In our notation, Keller showed that for $x_i(t)\in[0,1]$ evolving according to the equations
\[
x_i(t+1)= f(x_i(t))+\frac{\epsilon}{N}\sum_{j=1}^NA_{ij}f(x_j(t))
\] 
with $A_{ij}\in\{0,1\}$,
\[
f(x)=\alpha\left(\frac12-\left|x-\frac12\right|\right),
\]
and $\alpha$ satisfying certain assumptions ensuring uniform hyperbolic behavior, the system has an invariant absolutely continuous mixing measure. Again, the results are inferred using Markov partitions and coding. More on the study of GCMs in the context of ergodic theory can be found in Section \ref{Sec:ErgThCoupMaps}.

\medskip
In the next section we review some further results on the turbulent phase when $N\rightarrow \infty$, and a phenomenon that influenced the study of globally coupled maps.

\subsection{Violation of the law of large numbers}\label{Sec:ViolLN}
In \cite{kaneko1990globally} and \cite{kaneko1992mean}, Kaneko observed a phenomenon that he termed \emph{violation of the law of large numbers}  later also referred to as \emph{nonstatistical behavior} \cite{perez1992instabilities,sinha1992nonstatistical}. The observation is the following. Consider a system of globally coupled maps as in \eqref{Eq:CoupledLogisticMaps} with logistic local dynamics as in \eqref{Eq:LogisiticFam} and where the parameter $a$ is chosen so that for small enough $\epsilon$, the system exhibits turbulence/chaos. The mixing character of the local dynamics suggests that in the limit $N\rightarrow \infty$, the maps should be uncorrelated\footnote{More strongly $x_i(t)$ and $x_j(t)$ are expected to evolve independently of one another as if they were uncoupled.} and  the \emph{mean-field} coupling term 
\begin{equation}\label{Eq:Meanfield}
h_N(t):=\frac{1}{N}\sum_{j=1}^Nf(x_j(t))
\end{equation}
 should satisfy the law of large numbers and converge to a fixed value. For a finite system, one then would expect the time series $h_N(t)$ to be close to this fixed value plus some fluctuations going to zero as $N\rightarrow\infty$. To measure the size of the fluctuations, one can pick the Mean Square Deviation ($MSD$) defined as
\[
MSD(N):=\left\langle h^2_N -\langle h_N\rangle^2 \right\rangle
\]
where $\langle\cdot \rangle$ denotes the integral with respect to $P_N$ which is  the (unknown) distribution of $h_N$. In other words, $MSD(N)$ is the variance of $h_N$. In practice, $MSD(N)$ can be estimated from the time series $\{h_N(t)\}_{t\ge 0}$. If $\{x_i(t)\}_{i=1}^N$ were uncorrelated, one would expect $MSD(N)$ to decay as $N^{-1}$. Surprisingly, numerical simulations showed that after an initial decrease proportional to $N^{-1}$, for larger $N$ the quantity $MSD(N)$ stabilizes at a fixed small value ($10^{-1}$-$10^{-3}$) depending on the parameter $a$ -- see Figure 2 in \cite{kaneko1990globally}. These  fluctuations were deemed due to some ``coherence" among  units that would persist in the limit $N\rightarrow \infty$.  

Shortly after, in \cite{pikovsky1994globally}, a different point of view was put forward, and the observed lack of decay of $MSN(N)$ was imputed to the lack of stationarity in the system. The conclusion in \cite{pikovsky1994globally} was that in the limit $N\rightarrow \infty$ the system can be out of the equilibrium and wonder between different states on which the mean-field takes different values that account for the fluctuations observed in the time series of $h_N(t)$. More precisely, rather than describing the state of each map $x_i(t)$, one can investigate their distribution given by the measure
\begin{equation}\label{Eq:measuredist}
\mu_N(t):=\frac{1}{N}\sum_{j=1}^N\delta_{x_j(t)}
\end{equation}
and its evolution\footnote{In writing the evolution law we followed the notation in \cite{pikovsky1994globally}, but this is just another expression for the transfer operator of $f_{a(t)}$ applied to the measure $\mu(t)$
\begin{align*}
\mu_N(t+1)&=\left(f_{a_t}\right)_*\mu_N(t).
\end{align*}}
\begin{equation}\label{Eq:Selfcons1}
\mu_N(t+1)=\int \delta_{f_{a_t}(y)}\, d\mu_N(t)[y]
\end{equation}
where $f_{a_t}(y)$ is a map from a parametric family and
\begin{equation}\label{Eq:Selfcons2} 
a_t=a_0+\epsilon h_N(t),\quad\quad h_{N}(t)=\int y\,d\mu_N(t)(y)=\frac{1}{N}\sum_{j=1}^Nx_j(t)
\end{equation}
with $a_0\in \R$. One can consider equations \eqref{Eq:Selfcons1}-\eqref{Eq:Selfcons2} beyond the current setup substituting the empirical distribution $\mu_N(t)$ with any measure $\mu(t)$. In particular, if $\mu(t)=\lim_{N\rightarrow \infty}\mu_N(t)$\footnote{We will be more careful about the type of convergence here later on.}, one can think of equations \eqref{Eq:Selfcons1}-\eqref{Eq:Selfcons2} as describing the evolution of the system's state in the thermodynamic limit. In \cite{pikovsky1994globally}, evidence has been found that the dynamics of $\mu(t)$ is not necessarily asymptotic to a fixed point, i.e. the thermodynamic limit does not necessarily have a stable equilibrium state. Instead, one can imagine that  for $\mu(0)$ in some class of measures, the orbit  $\mu(t)$ can evolve towards, for example, a periodic orbit (with period $>1$), or even more complicated attractors. In this situation, one could expect that for $N$ sufficiently large, $\mu_N(t)$ would also be close to the attractor shadowing its dynamics and the orbit $\mu_N(t)$ would appear as a noisy version of an orbit on the attractor. Going back to the study of the mean-field $h_N(t)$, this suggests that rather than fluctuations around the expected value of the mean-field $\langle h_N\rangle$, one should consider fluctuations with respect to $\int f(y)\, d\mu(t)[y]$. An example in \cite{pikovsky1994globally} shows numerically that the violation of the law of large numbers can be resolved taking this point of view. The available examples of this kind usually arise when the local maps $f$ belong to a family of maps with nontrivial bifurcation structure. 

\medskip
Examples where \eqref{Eq:Selfcons1}-\eqref{Eq:Selfcons2} have a trivial attractor  given by an attracting fixed point  are also available. In fact,  if the local maps  have  uniformly hyperbolic properties -- e.g. they are smooth with uniform expansion -- and the coupling strength is small, the system is expected to have a unique equilibrium in the thermodynamic limit close to the SRB measure of one of the $f_{a}$, and $\mu(t)$ converges to this equilibrium exponentially fast, provided that the initial condition $\mu(0)$ is picked inside a suitable set of measures with some smoothness. The first instance where a claim of this kind has been rigorously proved is \cite{keller2000ergodic}. This was followed by many other results on the study of self-consistent operators which we are going  to review in Section \ref{Sec:StabFixedPoints}.

\subsection{Mean-field fluctuations and self-consistent transfer operators}
The observations of the violation of the law of large numbers, opened the way to the study of the nontrivial evolution of the mean-field $h_N(t)$ in the limit for $N\rightarrow \infty$. The main starting point for this study are the equations \eqref{Eq:Selfcons1}-\eqref{Eq:Selfcons2} that define a nonlinear self-consistent evolution law on the space of measures. The generator of this evolution law is also called  \emph{self-consistent operator} or \emph{nonlinear Perron-Frobenius operator}. This object was introduced in the setting of globally coupled maps first  in \cite{kaneko1992mean}, while another version was already used in \cite{kaneko1989towards} in the context of coupled map lattices. In continuous time, an analogous nonlinear evolution is known as a  \emph{nonlinear Fokker-Planck} (e.g. \cite{desai1978statistical}).

\medskip
In \cite{just1995bifurcations}, the author studied the self-consistent evolution in the case where the uncoupled map $f:[-1,1]\rightarrow[-1,1]$ belongs to the family of tent maps
\[
f_a(x)=1-a|x|
\] or to the logistic family \eqref{Eq:LogisiticFam}. He investigated: the linearization of  \eqref{Eq:Selfcons1}-\eqref{Eq:Selfcons2} around fixed points, the presence of periodic orbits for the evolution, and  formal conditions implying stability. Interestingly, the conditions are reminiscent of those appearing in the study of linear response for the uncoupled map $f$, suggesting that structural  and/or statistical stability are requirements for the existence of stable equilibria in the thermodynamic limit.
This is also supported by numerical evidence showing that in the logistic family, when the parameter $a$ is selected in an area where linear response for $f$ fails,  even a very small change in the coupling strength $\epsilon$ can produce notable effects in the observed dynamics. Similar considerations and further evidence have been also put forward in \cite{kaneko1995remarks} and \cite{chawanya1998bifurcation}. 

Nonstatistical behavior has been observed also in heterogeneous systems \cite{shibata1997heterogeneity}, i.e. systems where the local dynamics are not identical as in \eqref{Eq:GloballyCoupled}, but each map has different local dynamics.

In \cite{just1997collective}, the presence of stable periodic orbits for \eqref{Eq:Selfcons1}-\eqref{Eq:Selfcons2} is investigated for a system of coupled tent maps. A bifurcation digram in the parameters $(a,\epsilon)$ is obtained exhibiting a period-doubling cascade (see Figure 2 in \cite{just1997collective})\footnote{See also \cite{nakagawa1998collective}  for a thorough numerical study of the typical size of fluctuations in the mean-field across parameter space}.

\medskip
 Rather then on the self-consistent evolution of measures, some works focus only  on its effect on the evolution of the mean-field $h_N(t)$ -- for the relation between evolution of measure and $h_N(t)$ recall the second equality in \eqref{Eq:Selfcons2}. In \cite{ershov1995mean}, an analysis of the self-consistent equations is used to estimate that in a system of globally coupled tent maps, $h_N(t)$ has nontrivial fluctuations for certain values of the height of the tents and, most surprisingly, for any value of the coupling strength $\epsilon>0$. Numerical and analytic considerations estimate the fluctuations at $e^{-C\epsilon^{-2}}$. A similar analysis is carried out for logistic maps in \cite{ershov1997mean} estimating the size of the fluctuations at the much larger order of magnitued $O(\epsilon)$.

\smallskip
Some insight on the origin of the fluctuations can be obtained from the return plots depicting $h_N(t+1)$ versus $h_N(t)$ -- see \cite{chawanya1998bifurcation} for coupled tent maps and \cite{shibata1998tongue} for coupled logistic maps. These present a variety of characteristics depending on the local maps and strength of interactions. For example, the points $(h_N(t),h_N(t+1))$:
\begin{enumerate}
\item[i.] can be concentrated on a finite collection of points, see Figure 1b in \cite{chawanya1998bifurcation}, that can occur e.g. when \eqref{Eq:Selfcons1}-\eqref{Eq:Selfcons2} have an attracting periodic orbit;
\item[ii.] can lay close to a one-dimensional curve, like a circle, in which case $\{h_N(t)\}_{t\ge 0}$ shows quasi-periodic behavior, see Figure 1a in \cite{chawanya1998bifurcation} and Figure 1b in \cite{shibata1998tongue};
 \item[iii.] can present more complicated, but still low-dimensional structure, for example laying on what looks like the projection of a 2D torus in an higher dimensional space to the plane $(h_N(t),h_N(t+1))$, see Figure  1 (b) in \cite{nakagawa1999confined};
\item[iv.] or they can lack any evident structure whatsoever, see Figure 1a in \cite{shibata1998tongue}.
\end{enumerate}

\subsection{Lyapunov exponents for the mean-field dynamics}\label{Sec-LYap}

In order to characterize the patterns in i.-iv., several authors have put forward different approaches defining Lyapunov exponents associated to the time-series $\{h_N(t)\}_{t\ge 0}$.

\subsubsection{Lyapunov exponents of the self-consistent equations} 
In \cite{kaneko1995remarks}, Kaneko investigated the Lyapunov exponents of the self-consistent equations in \eqref{Eq:Selfcons1} for a system of globally coupled tent maps. Starting from a measure $\mu(0)$, a small perturbation is applied yielding $\mu'(0)=\mu(0)+\delta\nu$. Then the orbits $\mu(t)$ and $\mu'(t)$ are compared, and the top Lyapunov exponent is estimated for several perturbations $\delta\nu$. Situations as in point i. presented a negative exponent, confirming the presence of a periodic attractor for the self-consistent equations, while situations like iii. and iv. yielded a positive exponent. The estimated values of these exponents varying the parameters $(a,\epsilon)$ can be found in Figure 10 from \cite{kaneko1995remarks}.

\smallskip
Top Lyapunov exponents for the finitely many coupled tent maps and their self-consistent equations have also been studied in \cite{morita1997lyapunov}. 

\subsubsection{Collective Lyapunov exponents}

A different type of analysis has been proposed in \cite{shibata1998collective} where the focus is on the Lyapunov exponent of equations \eqref{Eq:CoupledLogisticMaps} when perturbing an initial condition along the direction of the mean-field \eqref{Eq:Meanfield} only. More precisely, considering an initial condition $(x_1(0),...,x_N(0))$, a perturbation of this initial condition along the diagonal direction is obtained putting
\[
x_i'(0)=x_i(0)+\frac{\delta}{N}
\]
for some $\delta\in[-\delta_0,\delta_0]$. One then obtains a Lyapunov exponent studying the rates of divergence (or convergence) of $h_N(t)$ and $h_N'(t)$ which are the mean-field along  the original and perturbed trajectory respectively. Crucially, the Lyapunov exponent estimated with this particular type of perturbation is independent of $N$, for $N$ sufficiently large, and can be much smaller than the value of the top  exponent for the whole system. This exponent is believed to detect information about the collective motion of the coupled system that is emergent in the thermodynamic limit. 

 \subsection{Ergodic theory of Globally Coupled Maps}\label{Sec:ErgThCoupMaps}
One of the first papers investigating the dynamic of globally coupled maps in the context of ergodic theory is \cite{keller1997mixing}, where conditions for existence and stability of a unique absolutely continuous invariant measure were established for a finite system of globally coupled tent maps.  In \cite{keller2000ergodic},  a similar result was proved for the thermodynamic limit with infinitely many globally coupled tent maps. Here the time evolution is given by a self-consistent transfer operator, and existence and uniqueness of a fixed measure was showed providing the first rigorous results on STOs. 

In \cite{jarvenpaa1997srb}, it was proved that $N$ analytic uniformly expanding weakly\footnote{I.e. when the interaction strength between maps is $\epsilon/N$ with $\epsilon$ small. }  coupled maps admit a unique a SRB measure $\mu_N$, independently of $N$, and $\mu_N$ converges to a limit $\mu$ for $N\rightarrow \infty$. 

 \cite{bardet2009stochastically} studies a STO undergoing a pitchfork bifurcation. Here the STO is of the type in \eqref{Eq:SelfconsistParamFam} with a family of piecewise uniformly expanding maps with two onto branches and  all maps preserving Lebesgue measure. It is shown that for small values of the coupling, Lebesgue is a stable fixed point for the STO. Increasing the coupling strength, Lebesgue looses stability, and two stable fixed measures appear. 

\cite{selley2016mean} studies existence and stability of fixed points for STOs  arising from systems of coupled doubling maps with piecewise linear interactions for different regimes of coupling strength. These results were generalized in \cite{balint2018synchronization} to a wider class of uniformly expanding maps. In \cite{selley2021linear} the authors studied linear response of fixed points for smooth uniformly expanding maps with smooth interactions. In \cite{galatolo2022self}, a general functional analytic framework to study fixed points of self-consistent transfer operators and their stability is provided. In \cite{bahsoun2022globally}, STOs arising from coupled Anosov diffeomorphisms have been considered. A more careful discussion of the above results is given in Section \ref{Sec:STOs} where we focus on self-consistent transfer operators. 

Most of the existence results in the papers above require small coupling strength. Increasing the coupling is expected to destroy the stability of the fixed points, and to eventually lead to clustering and synchronization for larger coupling strength. Rigorous studies of the bifurcations happening in between are unavailable. 

Before the onset of clustering and the related decrease in dimensionality of the attractors, the system can undergo a bifurcation via \emph{breaking of ergodicity} where multiple attractors of full dimension form, i.e from a situation where the finite dimensional system has a unique  absolutely continuous invariant probability (a.c.i.p.) measure, to a situation where the system has multiple a.c.i.p. measures supported on disjoint sets of positive Lebesgue measure and full dimension.

Ergodicity breaking is often related to the breaking of some symmetries of the system. For example, in \cite{fernandez2014breaking} Fernandez considered a system of $N$ doubling maps coupled via piecewise affine diffusive interactions and, for $N=3$, provided numerical evidence and rigorous arguments showing that, by increasing the coupling strength, the unique a.c.i.p. measure of the system breaks into multiple asymmetric ergodic measures having support of positive Lebesgue measure. It is important to notice that this happens when the system of coupled maps is still uniformly expanding. The discontinuities in the coupling are therefore to be considered responsible for the bifurcation. Breaking of ergodicity in a system of 3 globally coupled maps has been rigorously studied also in \cite{selley2016mean}, and in a system of 4 coupled maps has been studied in \cite{selley2018symmetry}.  \cite{fernandez2020computer}, \cite{fernandez2022conditioning} study algorithms to obtain computer assisted proofs of breaking of ergodicity for piecewise affine uniformly expanding coupled maps in any dimension. The algorithms provide a way to check existence of forward invariant sets given by unions of polytopes that, given the expansivity assumptions, will be granted to support  a.c.i.p. measures. 


\subsection{Other systems of coupled maps}\label{Sec:Others}
In this section we briefly review other types of interacting systems in discrete time. We do not aim at completeness, but rather at highlighting some interesting aspects and pointing to some research trends and works that are relevant for the study of globally coupled maps. 
\subsubsection{Coupled Maps Lattices}\label{SeubSec:CML}
The system described by the equations \eqref{Eq:CML} is an example of CML. The main difference with globally coupled maps, is that there is a notion of distance among maps (e.g. given by a lattice structure) and interactions are local, i.e. they are only among nearby maps  or the interaction strength decays with the distance as, for example, in 
\[
x_i(t)=f(x_i(t))+\epsilon\sum_{j\in \Lambda} \psi(|i-j|) h(x_i(t),x_j(t))\quad\quad i\in\Lambda
\]  
where $\Lambda\subset \Z^d$, $|i-j|$ is the distance between nodes $i$ and $j$, and $\psi(|i-j|)\rightarrow 0$ when $|i-j|\rightarrow \infty$; typically 
\[
\psi(|i-j|)=e^{-\alpha|i-j|}\quad\quad\alpha>0
\]
or decays sufficiently fast so that  $\sum_{j\in \Lambda}\psi(|i-j|)$ is summable. 

Perhaps the most important difference with globally coupled maps is that even in the case of infinite $\Lambda$, each map feels a  nonzero -- $O(\epsilon)$ -- influence from some of the other maps, while in globally coupled maps the interaction strength among any two given units goes to zero when $N\rightarrow \infty$ and only the cumulative effect of many interactions has an effect on the dynamics.

\medskip
CMLs originated as discretized models of continuous spatially extended systems  such as fluids and systems of chemical reactions with diffusion. The book \cite{kaneko2001complex} reviews the behavior of CML as investigated in the physics literature. Most of what is reported below can be found there.

As for the study of globally coupled maps, the local dynamics mostly employed in the study of CMLs are logistic and  tent maps to capture chaotic dynamics with intricate bifurcation structure.  Numerical studies  (\cite{kaneko1989pattern,kaneko1989spatiotemporal,grassberger1991phase,afraimovich1993travelling,gang1994spatiotemporally,afraimovich1994chaos,lambert1994stability,brindley1994spatio,chate1995spatiotemporal,kurths1995symmetry,bunimovich1996onset,chow1998traveling,franceschini2002quasiperiodic}) showed that CMLs exhibit a great variety of behaviors:
\begin{itemize}
\item[i.] \emph{Periodic behavior.} In this state, the lattice is divided in various connected domains grouping nearby sites. Within each domain, the sites have periodic dynamic with the same period. These states are  observed for example in coupled logistic maps on a 1D lattice with parameter in the doubling cascade window. The subdivision into domains depends on the initial condition. Different initial conditions lead to different domains with possibly different characteristic periods.  The number of possible domain configurations scales exponentially with the system size (this is analogous to the attractor crowding discussed in Section \ref{Sec:GlobCouplVarBehav}).
\item[ii.] \emph{Spatial bifurcations.} Starting from a configuration as described in i. and increasing the parameter of the logistic map slightly, one observes that the domains tend to remain intact, but the dynamics within each domain bifurcates. In particular, it first goes through a period doubling cascade, until it eventually becomes chaotic. Thus one ends up with orbits that are periodic on some domains and chaotic on others.

\item[iii.] \emph{Spatiotemporal Intermittency.} Further increase of the parameter can lead to destruction of the domains. In this case the dynamic looks non-stationary with each site alternating between stretches of time where it exhibits quasi-periodic behavior, and abrupt switches to erratic motion (temporal intermittency). Furthermore, at the same instant of time, some sites exhibit periodic behavior, while other show irregular dynamics (spatial intermittency).
\item[iv.] \emph{Fully Developed Chaos.} Further increase of the parameter for the local maps makes every site undergo  chaotic  motion. The orbits at each site become uncorrelated on large scales.  In some cases transition to chaos happens for effect of the coupling. 
\item[v.] \emph{Travelling Waves.}\ These appear in the range of parameters for the logistic map discussed at points i. and ii., if the coupling strength is increased. In this case the domains are not invariant anymore, but they can move across space.
\end{itemize}


\medskip
CMLs with \emph{unidirectional coupling}, e.g. on a 1D lattice where each node receives an interaction only from its neighbors on the left, exhibit interesting phenomena not observed in systems with more general coupling \cite{schmuser2001non,keller2001unidirectionally}.

\medskip
Coupled map lattices have been extensively studied also using tools from ergodic theory. In this branch of the literature, the local dynamics are usually uniformly hyperbolic (e.g. uniformly expanding) and the coupling strength is weak.   \cite{chazottes2005dynamics} reviews early works in its introduction, and collects also several papers on the topic.

A seminal work is \cite{bunimovich1988spacetime}. Here a 1D lattice of coupled uniformly weakly expanding maps is considered. The evolution equations look like
\[
x_i(t+1)=(1-\epsilon)f(x_i(t))+\epsilon [f(x_{i-1}(t)+f(x_{i+1}(t))]
\]
where $f$ is a uniformly expanding map. It was expected that if there are only finitely many coupled maps (e.g. on a finite periodic 1D lattice) and $\epsilon$ was sufficiently small so that the resulting map was expanding, then the system had a unique absolutely continuous invariant measure. The question was if also the infinite system (e.g. on $\Z$), admitted a unique  SRB\footnote{In this infinite-dimensional setting an SRB measure is, roughly speaking, a measure whose projection on finite dimensional marginals corresponding to finite subsets of the lattice are absolutely continuous. For a discussion on the definition of SRB measures for infinite-dimensional coupled maps see \cite{jarvenpaa2001definition}} measure for $\epsilon$ small but different from zero. In \cite{bunimovich1988spacetime} this question was answered in the affirmative with the use of symbolic dynamics.

Lattices of uniformly expanding coupled maps and their invariant measures have been later on studied in: \cite{gundlach1993spatio,gundlach1993spatioII,gundlach1993spatioIII,volevich1994construction,jarvenpaa1999note,schmitt2004spectral};  \cite{keller1992transfer} for piecewise expanding systems; \cite{bricmont1995coupled}, for analytic maps;   \cite{bricmont1996high}, using a high temperature expansion;  \cite{bricmont1997infinite}, that proved exponential decay of spatio-temporal correlations. 

SRB measures for coupled map lattices where the local dynamics has an hyperbolic attractor were studied using approaches from thermodynamics (e.g. polymer expansions) in \cite{pesin1988space, jiang1995equilibrium, jiang1998equilibrium, gundlach1993spatioII, bonetto2004analyticity} \cite{bricmont2013diffusion}.

Results on the spectral properties of Perron-Frobenius operators for various types of CML can be found in: \cite{keller1996coupled}; \cite{baladi1998spectrum}; \cite{fischer2000transfer} for analytic coupled maps using a cluster expansion;  \cite{rugh2002coupled}; \cite{jiang2003sinai} using the thermodynamic formalisms fo transfer operators;

Finally in \cite{keller2006uniqueness} a general framework for the study of Perron-Frobenius operators of coupled expanding maps has been put forward. In this paper the authors construct Banach spaces on the infinite dimensional phase space and a direct proof of the presence of a spectral gap is provided. The argument exploits the uniform expansion of the uncoupled dynamics, and the local nature of the interactions.

In \cite{maes1997stochastic} stochastic stability of the Gibbs states is investigated, while \cite{jiang2000smooth} studies linear response.  

Some works study finer statistical properties of CMLs: \cite{bardet2002limit}  investigates limit theorems; \cite{bardet2004spatio} studies large deviations;  escape rates in coupled map lattices with holes   are studied in \cite{bardet2011extensive} using symbolic dynamics, and in  \cite{faranda2018extreme} using the perturbation theory of transfer operators (with applications to synchronization). 

Entropy is  studied in \cite{dolgopyat1997entropy}. 

Increasing the coupling strength, the picture with a unique SRB measure is destroyed and one witness the appearance of: multiple Gibbs states \cite{keller1992some,losson1995phase,blank1997generalized,gielis2000coupled,bardet2006phase}; coherent structures \cite{bunimovich1990emergence,blank2013emergence}. See \cite{bunimovich1992coupled} for an example with simple uncoupled dynamics and coupling, where a full picture concerning bifurcations is rigorously established. Phase transitions and bifurcations in CML are rigorously studied also in \cite{de2010phase}\cite{bunimovich1998localized}. 

\medskip
\cite{just1998analytical,afraimovich2000topological,just2001equilibrium} focus on the topological properties of piecewise affine CMLs (rather than the measure theoretic ones presented above) using symbolic dynamics. 

\medskip In \cite{keller2009map} and \cite{bahsoun2022map}, the authors study maps coupled by collisions. Here uniformly hyperbolic maps are coupled to each other by rare  but very strong interactions:  on most of the phase space the system is uncoupled apart from a small set where the interactions can be large.   

\medskip
Reviews  of CMLs can be found in \cite{bunimovich1991simple,bunimovich1997coupled,bunimovich2005coupled}.

\subsubsection{Coupled map networks}

Coupled maps with more general types of coupling structures have been considered in the literature. Usually the maps are assumed to occupy the vertices of a graph and the presence of an edge prescribes an interaction. These systems have been considered in \cite{koiller2010coupled} and were termed \emph{coupled map networks}. Here the coupled maps are smooth and uniformly expanding, and the interactions (among maps connected by an edge) are piecewise affine. The equations describing the evolution of $N$ coupled maps on the 1D torus can be written as
\[
x_i(t+1)=f(x_i(t))+\epsilon\sum_{j=1}^NM_{ij}[x_j(t)-x_i(t)]\mod 1\quad\quad i=1,...,N
\] 
where $M_{ij}$ is a matrix of weights associated to each directed edge from node $j$ to node $i$. The paper provides sufficient conditions involving the matrix $(M_{ij})$  for the resulting dynamics to be piecewise hyperbolic.

\medskip
 In  \cite{pereira2020heterogeneously} the authors study uniformly expanding coupled maps on heterogeneous networks with evolution equations 
  \[
 x_i(t+1)=f(x_i(t))+\frac{1}{\Delta}\sum_{j=1}^NA_{ij}h(x_i(t),x_j(t))\mod 1\quad\quad i=1,...,N
 \]
where $(A_{ij})$ is the adjacency matrix of an heterogeneous graph, i.e. having most of the nodes making very few connections (low degree nodes) and a few nodes (called hubs) being connected to a large number of nodes\footnote{This graphs are akin to scale-free networks \cite{barabasi2009scale}.}. The parameter $\Delta$ is the maximum in-degree of the network. Here it was showed that a mean-field reduction can be made for the dynamics of the hub nodes where the average of the interaction is substituted by an expectation, and the reduction holds for times exponentially large in the system's size.

\medskip
For another treatment of the effect of the structure of interactions in shaping the dynamics see \cite{afraimovich2010dynamical}.


\section{Self-Consistent Transfer Operators}\label{Sec:STOs}

 As we have argued above, self-consistent transfer operators (STOs) arise as thermodynamic limit of coupled maps. More generally, to define a self-consistent operator acting on measures one needs to specify a mapping that to each measure associates a linear operator on measures; the STO than acts  taking a measure, associating the corresponding linear operator to it, and applying this operator to the measure itself. This is made precise in the following definition.
\begin{definition}\label{Def:Selconsistentoperator}
Given a Borel space $(X,\mc B)$\footnote{In what follows,  $X$ is going to be a compact metric space.}, let's denote by $ \mc M(X)$ the set of finite signed Borel measures on $X$ and $V\subset\mc M(X)$ a subspace.  

Denote by $\End_p(V)$ the set of linear endomorphisms of $V$ preserving the total measure of $X$, i.e. such that for every $A\in\End_p(V)$ and $\mu\in V$, \[A[\mu](X)=\mu(X).\] 

A mapping $T:V \rightarrow \End_p(V)$ defines the \emph{self-consistent operator} $\mc T:V\rightarrow V$  as
\[
\mc T(\mu):=T(\mu)[\mu]
\]
with the above notation standing for the operator $T(\mu)$ applied to the measure $\mu$.
\end{definition}
 $\mc T$ is a nonlinear selfmap of $V$. Depending on the context, the object defined above has been given different names.  In the context of nonlinear Markov chains it is  referred to as \emph{nonlinear Perron-Forbenius operator}.  
 
 The main goal  is to study the properties of $\mc T$  from knowledge of the mapping $T$. 
Notice that given any map $\mc P:\mc M(X)\rightarrow \mc M(X)$ such that $\mc P(\mu)(X)=\mu(X)$ for all $\mu\in \mc M(X)$, there exist (many)  mappings $T:\mc M(X)\rightarrow \End_p(\mc M(X))$ such that the associated self-consistent operator $\mc T$ equals $\mc P$. It is therefore crucial to restrict  to some specific classes of $T$  to obtain self-consistent operators amenable to study. Below we list some possible setups. 

\begin{itemize}
\item {\bf Average of self-consistent operators.}
 Consider a measurable map $\gamma:X\rightarrow \End_p(V)$, and define $T:V\rightarrow \End_p(V)$ as
 \[
T(\mu):=\int_{X}\gamma(x)\,d\mu(x)
 \] 
 and 
 \[
 \mc T\mu=T(\mu)[\mu]=\int_Xt(x)[\mu]\,d\mu(x).
 \]
 The above can be interpreted as an average of transfer operators, where the average is with respect to the measure it's applied to.
\item {\bf Nonlinear Markov Chains.}  As a particular example of the above, consider $P:\mc B\times X \times X\rightarrow \R^+_0$ such that for every $x,y\in X$, $P(\cdot,x,y):\mc B\rightarrow \R^+_0$ is a probability measure and for all $P(B,\cdot,\cdot)$ is measurable for all $B\in \mc B$. $P$ should be interpreted as a $y$ dependent transition probability. Define $T: \mc M(X)\rightarrow \End_p(\mc M(X))$ as 
\[
T(\mu)[\nu](\cdot)=\int_X\left(\int_XP(\cdot,x,y)\,d\mu(y)\right)d\nu(x)
\]
and
\[
\mc T\mu(\cdot)=\int_X\int_XP(\cdot,x,y)\,d\mu(y)d\mu(x).
\]
\item {\bf Globally mean-field coupled maps with full permutation symmetry. } For simplicity let $X=\T=\R\backslash\Z$, and consider functions $f:\T\rightarrow \T$ and $h:\T\times\T\rightarrow \R$ and  the system of globally coupled maps given by
\begin{equation}\label{Eq:EqsCoupledSyst}
x_i(t+1)=f(x_i(t))+\frac{1}{N}\sum_{j=1}^Nh(x_i(t),x_j(t))\mod 1\quad\quad i=1,...,N
\end{equation}
where $x_i(t)$ describes the state at time $t$ of the $i$-th unit. Defining for any $\mu\in\mc M(\T)$, $f_{\mu}:\T\rightarrow\T$ as   
\[
f_{\mu}(x)=f(x)+\int_{\T}h(x,y)\,d\mu(y) \mod 1
\]and $\mu^{(N)}_t:=\frac{1}{N}\sum_{j=1}^N\delta_{x_i(t)}$ one has
\[
\mu^{(N)}_{t+1}=\left(f_{\mu^{(N)}_t}\right)_*[\mu^{(N)}_t]
\]
where $(f_{\mu^{(N)}_t})_*$ denotes the push-forward of $f_{\mu^{(N)}_t}$.  This leads  to the definition of $T:\mc M(\T)\rightarrow\End_p(\mc M(\T))$ as  $T(\mu)[\nu]=(f_{\mu})_*\nu$ and 
\begin{equation}\label{Eq:STOThermLimit}
\mc T\mu=(f_{\mu})_*[\mu].
\end{equation}
 Under some continuity assumptions on $f$ and $h$, one can see that if $\mu^{(N)}_t$ converges weakly to $\mu$ when $N\rightarrow \infty$, then $\mu^{(N)}_{t+1}$ converges weakly to $\mc T\mu$ (see e.g. the introduction of \cite{selley2021linear} for a discussion). In this sense, $\mc T$ describes the thermodynamic limit of the system.
\item {\bf Globally coupled maps without  symmetry.} Again let $X=\T=\R\backslash\Z$, given functions $f_i:\T\rightarrow \T$ and $h_{ij}:\T\times\T\rightarrow \R$,  consider the system of globally coupled maps given by
\[
x_i(t+1)=f_i(x_i(t))+\frac{1}{N}\sum_{j=1}^Nh_{ij}(x_i(t),x_j(t))\mod 1\quad\quad\forall i=1,...N.
\]
For any $\mu\in \mc M(\T^N)$ let $F_{\mu,i}:\T\rightarrow \T$ be
\[
F_{\mu,i}(x)=f_i(x)+\int_{\T^N}\left(\frac{1}{N}\sum_{j=1}^Nh_{ij}(x,y_j)\right)\,d\mu(y_1,..,y_j)\mod 1
\]
and $\bo F_\mu:\T^N\rightarrow \T^N$ given by $\bo F_\mu=(F_{\mu,1},...,F_{\mu,N})$. Then define $T:\mc M(\T^N)\rightarrow \End_p(\mc M(\T^N))$ as
\begin{equation}\label{Eq:SelfConsistNoSymm}
T(\mu)[\nu]:=(\bo F_\mu)_*\nu.
\end{equation}
The corresponding $\mc T$ is an extension of the one at the point above and becomes the previous one in the case where $f_i=f$ and $h_{ij}=h$ for all $i,j\in[1,N]$, and $\mu=\mu_1\otimes...\otimes\mu_1$ is a product measure with all identical factors.
\item {\bf Parametric families of maps.} Given a parametric family of maps on $X$, $\{f_\gamma\}_{\gamma\in \Gamma}$,  and  $\hat\gamma:\mc M(X)\rightarrow \Gamma$, then one can define
$T(\mu)[\nu]=(f_{\hat\gamma(\mu)})_*\nu$
and \begin{equation}\label{Eq:SelfconsistParamFam} \mc T\mu=(f_{\hat\gamma(\mu)})_*\mu.\end{equation} 
 \end{itemize}

Below, we are going to illustrate some of the main available techniques for the analysis of self-consistent transfer operators and their stable fixed points (Section \ref{Sec:StabFixedPoints}). We will then discuss numerical and rigorous results on linear response for fixed points of STOs and globally coupled maps (Section \ref{Sec:LinResp}),  and some further directions in the study of STOs when their attractors are different from fixed points (Section \ref{Sec:TwardsmreComp}).

\subsection{Stability for fixed points of STOs}\label{Sec:StabFixedPoints}

In this section we describe the main frameworks used to study stability and convergence to fixed points of STOs arising from globally coupled maps. The objective is not to give the results in their most general formulations, but restrict to a simple example where only the core ideas of each reviewed framework are highlighted. 

The example is the following: Consider a system of coupled maps with $x_i(t)\in\T$ and
\[
x_i(t+1)=f\left(x_i(t)+\frac{\epsilon}{N}\sum_{j=1}^Nh(x_i(t),x_j(t))\mod 1\right)\quad\quad i=1,...,N
\]where  $f(x)=2x\mod1$ (the doubling map), and $h$ is some smooth coupling function. The corresponding STO is
\begin{equation}\label{Eq:STO}
\mc T_\epsilon \mu = \mc P\mc L_{\epsilon,\mu}\mu
\end{equation}
with $\mc P$  the  linear transfer operator of the doubling map, which on $L^1(\T)$ acts as
\begin{equation}\label{Eq:TransfOpDoubMap}
\mc P\phi(x)=\frac{1}{2}\phi\left(\frac x2\right)+\frac12\phi\left(\frac{x+1}{2}\right),
\end{equation} and $\mc L_{\epsilon,\mu}$  the transfer operator associated to the mean-field coupling, i.e.  the transfer operator of the map
\[
g_\mu(x):= x+\epsilon \int h(x,y)d\mu(y)\mod 1.
\]

We are going to review three methods to study fixed points of the family of STOs above. The first two are devised for the case of small coupling, while the last one treats some situations that can arise in the case of strong coupling. These are: a \emph{functional analytic approach} that extends the spectral gap properties of linear operators to STOs (Section \ref{Sec:FuncAnalyticApp}); an approach with \emph{convex cones} that studies the contraction properties of STOs with respect to the Hilbert projective metrics (Section \ref{Sec:ConeApproach}); and an approach devised to study  \emph{synchronized states} (Section \ref{Sec:SyncedStates}).

\subsubsection{Functional analytic approach}\label{Sec:FuncAnalyticApp}
This is the approach that under different forms was used e.g. in \cite{keller2000ergodic,selley2016mean,balint2018synchronization}. 

The strategy consists of the following steps:

\smallskip
{\it Step 1.} Use Schauder's fixed-point theorem to prove that $\mc T_\epsilon$ has a fixed point $\mu_*$.

{\it Step 2.} Use the spectral gap of the family of linear operators $\mc P\mc L_{\epsilon,\mu}$ and continuity of $\mu\mapsto \mc P\mc L_{\epsilon,\mu}$  to show that the fixed point is attracting when $\epsilon$ is sufficiently small.

\medskip
Since  we are only going to deal with absolutely continuous measures, if $\mu$ has density $\phi$ we use notations: $\mc L_{\epsilon,\phi}$ and $g_\phi$.

{\it Step 1.} Consider the set 
\[
B_{L}:=\left\{\phi:\T\rightarrow \R^+_0:\, |\phi|_{Lip}\le L,\,\int_\T\phi(x)dx=1\right\}
\]
where $|\phi|_{Lip}=\sup_{x\neq y}\frac{\phi(x)-\phi(y)|}{|x-y|}$ denotes the Lipschitz semi-norm. The first thing to notice is that for $\epsilon>0$ sufficiently small, there is $L$ for which $B_L$ is forward invariant under action of $\mc T_\epsilon$.

\begin{lemma}\label{Lem:InvarianceBallLip} There is $\epsilon_0>0$ such that for any $|\epsilon|<\epsilon_0$ there is $L_0=O(\epsilon)$, such that for every $L>L_0$
\[
\mc  T_\epsilon B_L\subset B_L.
\]
\end{lemma}
\begin{proof}
See Appendix \ref{App:Comput}.
\end{proof}
Furthermore $\mc T_\epsilon$ is continuous.
\begin{lemma}\label{Lem:ContC0Top} With the parameters as in Lemma \ref{Lem:InvarianceBallLip}, $\mc T_\epsilon$ is 
is continuous in the $C^0$ topology\footnote{By $C^i$ topology we mean the topology generated by the norm
\[
\|u\|_{C^i}:=\sum_{j=0}^i\sup\left|\frac{d^j}{dx^j}u\right|
\]
for $u\in C^i$.}. 
\end{lemma}
\begin{proof}
See Appendix \ref{App:Comput}.
\end{proof}
Since $B_L$ is a convex, compact (in $C^0$) set, and $\mc T_\epsilon$ is continuous, by Schauder's fixed point  theorem $\mc T_\epsilon$ has a fixed point $\phi_*\in B_L$. 

\medskip
{\it Step 2.} This step is a bit more involved. First of all, one needs to modify Step 1. a bit to find a forward invariant set of functions more regular than just Lipschitz, for example $C^2$ with uniformly bounded first and second derivative:
\[
\mc B_{L_1,L_2}:=\left\{\phi\in C^2(\T,\R):\,|u'|\le L_1,\,|u''|\le L_2,\,\int \phi=1\right\}.
\]
In particular, arguments as those presented in Step 1 allow to conclude that the set $\mc B_{L_1,L_2}$ is forward invariant for suitable values of $L_1$ and $L_2$, and $\phi_*$ is in fact $C^1$ and has Lipschitz derivative with Lipschitz constant bounded by $L_2$.

Then, for every function $u=\phi_1-\phi_2$ where $\phi_1,\phi_2\in \mc B_{L_1,L_2}$  one proves that 
\begin{align}
\| \mc P \mc L_{\epsilon,\phi_*} u \|_{C_1}&\le \alpha\|u\|_{C_1}\label{Eq:Contraction}\\
\| \mc P(\mc L_{\epsilon,\phi} -\mc L_{\epsilon,\phi_*}) \phi_*\|_{C_1}&\le K \epsilon\cdot \|\phi-\phi_*\|_{C_1} \label{Eq:Deponphi}
\end{align}
with $\alpha\in(0,1)$, $K\ge 0$ depending on $L_1$ and $L_2$. Equation \eqref{Eq:Contraction} is a spectral gap condition for the linear operator $ \mc P \mc L_{\epsilon,\phi_*}$ which is implied (for $|\epsilon|$ sufficiently small) in a standard way by the uniform expansion and smoothness of the map $f\circ g_{\phi_*}$. Equation \eqref{Eq:Deponphi} is a continuity relation for the family of operators $\{\mc L_{\epsilon,\phi}\}$ and is proven in Lemma \ref{Lem:Continuityrelation}. 
\begin{lemma}\label{Lem:Continuityrelation}
Assume that $h$ is  $C^2$ and $\phi_*$ is as above, then there is $K\ge 0$ such that for all $\phi\in \mc B_{L_1,L_2}$
\[
\|\mc P( \mc L_{\epsilon,\phi} -\mc L_{\epsilon,\phi_*}) \phi_*\|_{C_1}\le K \epsilon\cdot \|\phi-\phi_*\|_{C_1}.
\]
\end{lemma}
\begin{proof}
See Appendix \ref{App:Comput}.
\end{proof}


By triangular inequality
\begin{align*}
\|\mc P \mc L_{\epsilon,\phi}(\phi)-\mc P \mc L_{\epsilon,\phi_*}(\phi_*)\|_{C^1} 
&\le \|\mc P \mc L_{\epsilon,\phi_*}(\phi-\phi_*)\|_{C^1}+\|\mc P (\mc L_{\epsilon,\phi_*}-\mc L_{\epsilon,\phi})\phi_*\|_{C^1}.
\end{align*}
The first term is less than $\alpha\|\phi-\phi_*\|_{C^1}$ by Eq. \eqref{Eq:Contraction};  the second term is bounded by $ O(\epsilon)\|\phi_*-\phi\|_{C^1}$ by Eq. \eqref{Eq:Deponphi}.
Putting the above estimates together
\begin{equation}
\|\mc P \mc L_{\epsilon,\phi}(\phi)-\phi_*\|_{C^1} \le \left[\alpha+ O(\epsilon)\right]\|\phi-\phi_*\|_{C^1}
\end{equation}
and for $|\epsilon|$ sufficiently small, 
\[
\alpha+ O(\epsilon)<1
\]
and one gets the desired contraction which implies that $\mc T_\epsilon^n\phi \rightarrow \phi_*$.

In the above setup one can  also obtain estimates of 
$
\|\phi_{\epsilon*}-\phi_{\epsilon'*}\|_{C^1}
$
where $\phi_{\epsilon*},\,\phi_{\epsilon'*}$ are the fixed points for $\mc T_{\epsilon}$ and $\mc T_{\epsilon'}$. In Section \ref{Sec:LinResp} we discuss some approaches to obtain differentiability of the mapping $\epsilon\mapsto\phi_{\epsilon*}$.

\subsubsection{The cone approach}\label{Sec:ConeApproach} In this approach, rather than studying the STO with respect to the norm of some Banach space, one considers its action with respect to the Hilbert projective metric on some convex cones of functions (see e.g. \cite{liverani1995decay}). To study the STO in \eqref{Eq:STO}, we can restrict its action to the cone of $\log$-Lipschitz functions 
\[
\mc V_a:=\left\{\phi:\T\rightarrow \R^+:\,\frac{\phi(x)}{\phi(y)}\le e^{a|x-y|}\right\}
\]
for some $a>0$, which is endowed with the Hilbert projective metric $\theta_a$. The peculiarity of the Hilbert metric is that any linear application between two convex cones is a contraction with respect to their Hilbert metrics. For this reason convex cones have been used to study the contraction properties of transfer operators of hyperbolic maps. 

It is immediate to check that 
\begin{equation}\label{Eq:ConEstDoub}
\mc P(\mc V_a)\subset \mc V_{a/2}.
\end{equation} For what concerns $\mc L_{\epsilon,\phi}$ we have the following
\begin{lemma}\label{Lem:COnes}
For all $\phi\in \mc V_a$, 
\begin{equation}\label{Eq:ConEstCoup}
\mc L_{\epsilon,\phi}\phi\in \mc V_{a'}
\end{equation}
with $a':=a[1+\mc O(\epsilon)]+\mc O(\epsilon)$.
\end{lemma}
\begin{proof}
See Appendix \ref{App:Comput}.
\end{proof}
Eq. \eqref{Eq:ConEstDoub} and Eq. \eqref{Eq:ConEstCoup} together imply that for suitable values of $a$ and $\epsilon$\footnote{More precisely if 
\[
\frac12\left[a(1+\mc O(\epsilon))+\mc O(\epsilon) \right]<a
\]
which can always be realized for $a$ sufficiently large.}, there is $\lambda\in[0,1)$ such that
\[
\mc T_\epsilon\mc V_a\subset \mc V_{\lambda a}.
\]
If $\mc T_\epsilon$ were linear, this would be enough to conclude that $\mc T_\epsilon$ is a contraction, and standard arguments would lead to existence of a fixed point together with uniqueness and stability results. However, $\mc T_\epsilon$ is nonlinear, so an extra step is needed to conclude the argument. In \cite{selley2021linear} we used the explicit expression for the Hilbert metric $\theta_a$ to prove that $\mc T_\epsilon:\mc V_a\rightarrow \mc V_{\lambda a}$ is a contraction when $|\epsilon|$ is sufficiently small. With some additional arguments one can use this fact to conclude that there exist a fixed density $\phi_*\in\mc V_a$ and that 
\[
\sup_{x\in \T}|\mc T_\epsilon^n\phi(x)-\phi_*(x)|\rightarrow 0
\]
exponentially fast.

\subsubsection{Synchronized states and the study of their stability}\label{Sec:SyncedStates}
Let's consider the explicit choice for interaction function
\[
h(x,y)=\sin(2\pi x)\cos(2\pi y).
\]
Notice that the measure $\delta_0$ is a fixed point for $\mc T_\epsilon$, in fact 
\[
g_{\delta_0}(x)=x-\epsilon \sin(2\pi x)
\] and  since 0 is a fixed point for $f$ and $g_{\delta_0}$, $\mc T\delta_0=f_*(g_{\delta_0})_*\delta_0=\delta_0$. 

The state $\delta_0$ can be seen as a synchronized state in the thermodynamic limit: Letting $(x_1,...,x_N)\in \T^N$ be the state of the  finite-dimensional system,  if
 \[
 \lim_{N\rightarrow \infty}\frac1N\sum_{i=1}^N\delta_{x_i}=\delta_0,
 \]
 when $N$ increases, the fraction of states $x_1,...,x_N$ that are further than any $\eta>0$ from zero must go to zero, i.e.
 \[
 \lim_{N\rightarrow \infty}\frac{\#\{x_i:\,|x_i-0|>\eta\}}{N}=0.
 \]
 
  Below we argue that if $\epsilon$ is in a certain range, then $\delta_0$ is stable in the following sense: There is $\Delta_0>0$ such that if $\mu\in \mc M_1(\T)$ has support contained in $[-\Delta_0,\Delta_0]$, then 
\[
\lim_{n\rightarrow \infty}\mc T_\epsilon^n\mu= \delta_0\quad\quad\mbox{weakly}.
\]
To show this, one can start by noticing that fixing $\epsilon$ in $(\frac1{2\pi},\frac3{2\pi})$, 0 is an attracting fixed point for the map $f_{\delta_0}$, and there are $\lambda'\in[0,1)$ and $I_{\delta_0}=(-\Delta', \Delta')$  such that $|f_{\delta_0}'(x)|<\lambda'$ for $x\in I_{\delta_0}$. By continuity, there is $\Delta_0>0$ sufficiently small and $\lambda\in[0,\lambda')$ such that for any  $\Delta<\Delta_0$  and $\mu$ with support contained in $[-\Delta,\Delta]$, $|f_\mu'|<\lambda$ on $[-\Delta,\Delta]$ and therefore the measure
\[
\mc T_\epsilon\mu=(f_\mu)_*\mu
\]
has support contained in $[-\lambda\Delta,\lambda\Delta]$. Arguing by induction, $\mc T_\epsilon^n\mu$ has support in $[-\lambda^n\Delta,\lambda^n\Delta]$, and for $n\rightarrow \infty$, $\mc L_\epsilon^n\mu$ converges weakly to $\delta_0$.

\medskip
In \cite{selley2022synchronization}, the picture above is generalized to the case where multiple clusters of coupled maps interact and the number of maps in each cluster goes to infinity. For example, taking a setup with two clusters, the states of the maps in each cluster are given by $x_1,...,x_N\in \T$ and $y_1,...,y_N\in \T$ and the evolution equations are
\begin{align*}
x_i(t+1)&=f^{(1)}\left( x_i(t)+\frac{1}{N}\sum_{j=1}^Nh_{11}(x_i(t),x_j(t))+\frac{1}{N}\sum_{j=1}^Nh_{12}(x_i(t),y_j(t))\right)\\
y_i(t+1)&=f^{(2)}\left( y_i(t)+\frac{1}{N}\sum_{j=1}^Nh_{22}(y_i(t),y_j(t))+\frac{1}{N}\sum_{j=1}^Nh_{21}(y_i(t),x_j(t))\right)
\end{align*}
where $f^{(1)}$ and $f^{(2)}$ are the local dynamics in the first and second cluster, while $h_{11}$, $h_{22}$, $h_{12}$, $h_{21}$ are respectively the interactions: among sites in the first cluster, among sites in the second cluster, from cluster 2 to cluster 1, and from cluster 1 to cluster 2.

The STO associated to the infinite limit of the system above is given by $\mc T:\mc M_1(\T^2)\rightarrow \mc M_1(\T^2)$ with
\[
\mc T\mu=(\bo F_\mu)_*\mu
\]
where $\bo F_\mu=(F_{\mu,1},F_{\mu,2}):\T^2\rightarrow \T^2$ is given by
\begin{align*}
F_{\mu,1}(x)&=f^{(1)}\left(x+\int_{\T^2}h_{11}(x,x')d\mu(x',y')+\int_{\T^2}h_{12}(x,y')d\mu(x',y')\right)\\
F_{\mu,2}(y)&=f^{(2)}\left(y+\int_{\T^2}h_{22}(y,y')d\mu(x',y')+\int_{\T^2}h_{21}(y,x')d\mu(x',y')\right).
\end{align*}
Notice that when $\mu=(\delta_x,\delta_y)$, then
\[
\mc T(\delta_x,\delta_y)=(\delta_{F_{(\delta_x,\delta_y),1}(x)},\delta_{F_{(\delta_x,\delta_y),2}(y)})
\]
and therefore the map  $G:\T^2\rightarrow\T^2$ 
\begin{equation}\label{Eq:G}
G(x,y)=(F_{(\delta_x,\delta_y),1}(x),F_{(\delta_x,\delta_y),2}(y)),
\end{equation}
prescribes the evolution of $(\delta_x,\delta_y)$. In \cite{selley2022synchronization},  sufficient conditions involving $G$ are given for the STO $\mc T$ to have stable fixed synchronized states. That paper considers also setups with multiple clusters where the STO has a stable fixed point which is a product of delta and absolutely continuous measures: this means that some clusters are in a synchronized state, while others are in a turbulent state. The clusters can be chosen so that the equations describing the system have full symmetry giving rise to what is sometimes called a chimera state (see Section \ref{subsec:chimeras}). 

\subsection{Linear response}\label{Sec:LinResp}  Given a high-dimensional system composed of many interacting units, how does its behavior change if the dynamics of its components is perturbed slightly? In particular, if the dynamics of the components is perturbed with a perturbation of magnitude $\epsilon$, is the change in the global behavior of the system still of order $\epsilon$? If yes the system is said to have linear response.

Linear response of high-dimensional systems has important relations to the study of climate models \cite{lucarini2018revising}. It is generally believed that high-dimensional systems exhibit linear response of their physical relevant measures. This is in conjunction with the chaotic hypothesis of Gallavotti and Cohen \cite{gallavotti1995dynamical} stating that high-dimensional systems are akin to Axiom A systems, for which linear response is known to hold \cite{ruelle2009review}. In the works we review below the question of linear response is addressed in some setups of globally coupled maps and STOs.

\subsubsection{Linear response for attracting fixed points of STOs}
Some of the  works mentioned above  give sufficient conditions for the  stable fixed point  $\phi_{\epsilon*}$ of  a parametric family of STOs $\mc T_\epsilon$ to be differentiable in $\epsilon$, and provide a linear response formula. 

In \cite{selley2021linear}, the cone approach yields differentiability of $\epsilon\mapsto \phi_{\epsilon*}$ from an interval $(-\epsilon_0,\epsilon_0)$ to $C^1$ densities. The main idea is to consider curves $\gamma:(-\epsilon_0,\epsilon_0)\rightarrow C^k(\T,\R)$ -- for a sufficiently large $k$ -- and the action $\mc T$ on these curves  given by
\[
(\mc T\gamma)(\epsilon):=\mc T_\epsilon(\gamma(\epsilon)),\quad\quad\forall \epsilon\in(-\epsilon_0,\epsilon_0).
\]
Loosely speaking,  the strategy consists in restricting to an invariant class of curves for   the action  $\mc T$, and use  Schauder's fixed point theorem to  show existence of an invariant curve $\gamma_*$ for $\mc T$ with the sought after differentiability property and such that $\gamma_*(\epsilon)$ is in the invariant cone for $\mc L_\epsilon$. By the discussion in Sect. \ref{Sec:ConeApproach}, this fixed curve $\gamma_*$ must satisfy
\[
\gamma_*(\epsilon)=\phi_{\epsilon *}
\]
and therefore $\phi_{\epsilon*}$ has a differentiable dependence on $\epsilon$. Once differentiability has been established, one can exhibit a linear response formula for the derivative of $\phi_{\epsilon*}$ with respect to $\epsilon$.

In \cite{galatolo2022self}, sufficient conditions are given in terms of the spectral properties of the linear operator   for the uncoupled system, and of the derivative of the nonlinear family of STOs $\{\mc T_\epsilon\}_\epsilon$ with respect to $\epsilon$. The main requirements are that:

i) $\mc T_\epsilon:\mc B_i\rightarrow \mc B_i$ for $i\in\{w,s,ss\}$ corresponding to three Banach spaces $\mc B_{w}\supset\mc B_{s}\supset \mc B_{ss}$ with norms $\|\cdot\|_w\le \|\cdot \|_s\le \|\cdot \|_{ss}$;

ii) the resolvent of the linear operator $\mc P$, $(\Id-\mc P)^{-1}$ is bounded on densities with zero integral from $\mc B_w$;

iii) an assumption that loosely speaking requires that $\mc T_\epsilon$ is Lipschitz in $\epsilon$\footnote{As an application from a subset of $\mc B_s$ to $\mc B_w$.} and differentiable in $\epsilon$ at $\epsilon=0$.

Under these assumptions $\epsilon\mapsto\phi_{\epsilon*}$ is differentiable at zero with respect to the weak norm $\|\cdot\|_w$, more precisely
\[
\lim_{\epsilon\rightarrow 0}\left\|\frac{\phi_{\epsilon*}-\phi_{0*}}{\epsilon}-(\Id-P)^{-1}\,\left.\frac{d}{d\epsilon}\mc L_\epsilon\right|_{\epsilon=0}\,\phi_{0*}\right\|_w=0.
\]

\subsubsection{Linear response for heterogeneous systems}\label{Sec:LinRespHet}
The setup and results below can be found in \cite{wormell2018validity,wormell2019linear}. Here linear response is studied for  systems of  different (finitely or infinitely many)  interacting units  belonging to a family of maps that does not necessarily exhibit linear response: the state of the global system at time $t$ is $\bo x(t)=(x_1(t),...,x_N(t))\in M^N$, and  the time evolution is given by
\begin{equation}\label{Eq:MicWormellSetup}
x_i(t+1)=f\left(x_i(t);\, a_i;\, \Phi(\bo x(t));\, \epsilon\right)\quad\quad i=1,...,N
\end{equation}
with 
\begin{equation}\label{Eq:MapsCoupledMeanfield}
f(x;a_i,\Phi,\epsilon)=f_{a_i}(x)+h_{\Phi}(x)+\epsilon g(x)\quad\quad a_i,\,\Phi,\,\epsilon\in \R.
\end{equation}
 The maps $f_{a_i}$ are a version of the logistic map with parameter $a_i$\footnote{More precisely is a skew-product, $f_{a_i}:\T\times\T\rightarrow\T\times\T$ with base the doubling map and fiber maps being either the identity or the logistic map with parameter $a_i$.}; $h_{\Phi(\bo x)}$ is a mean-field interaction term where
  \begin{equation}\label{Eq:MeanfieldpAram}
 \Phi(\bo x(t))=\frac{1}{N}\sum_{i=1}^N\phi(x_i(t))
 \end{equation}
is a mean-field parameter; and $\epsilon g(x)$ is a perturbation depending on the parameter $\epsilon$ w.r.t. which linear response is investigated. Different units can have different values of $a_i$ which are assumed to be drawn i.i.d. with respect to some (smooth enough) probability distribution. Crucially, depending on the distribution of the parameters $a_i$, the physical measure(s) of the maps $f_{a_i}$ may present or fail to exhibit linear response.

Considering a global observable\footnote{Global observables depend on the states of all the components, but weakly. Averages over the states of the components are instances of global observables. The temperature of a gas, for example, is a global observable  proportional to the average of the kinetic energy of the particles constituting the gas.}
\[
\Psi_N(\bo x):=\frac{1}{N}\sum_{i=1}^N\psi(x_i),
\]
and letting $\mu_\epsilon$ on $M^N$ be a physical invariant measure of the system, one wonders whether $\epsilon\mapsto \mb E_{\mu_\epsilon}[\Psi_N]$ is differentiable. 

For example, in the case where there is no mean-field coupling, i.e. $h_{\Phi}=0$, then
\[
\mb E_{\mu_\epsilon}[\Psi_N]=\frac{1}{N}\sum_{i=1}^N\mb E_{\mu_{\epsilon,i}}[\psi]
\]
where $\mu_{\epsilon,i}$ are the physical measures for the maps $f_{a_i}(x)+\epsilon g(x)$. It is immediate that if the maps  $f_{a_i}$ all satisfy linear response when perturbed adding $\epsilon g(x)$ to their equations, then so will the global uncoupled system. More surprisingly, it is shown in \cite{wormell2018validity} that even if the local dynamics don't satisfy linear response, the global finite-dimensional system  does, provided that the distribution of the $a_i$ is sufficiently smooth. In the same paper, it is shown that for some singular distributions of the $a_i$, linear response fails for the global system. 


Another interesting finding in \cite{wormell2019linear} is that when there is a mean-field coupling among the units, i.e. $h_\Phi\neq 0$, the authors bring numerical and analytical evidence showing that even if the microscopic units  exhibit linear response, in the thermodynamic limit, the global system can fail to do so. This can be rephrased by saying that even if the microscopic components leading to the definition of a STO as in \eqref{Eq:STOThermLimit} exhibit linear response, the dynamic of the STO might not satisfy linear response. Examples are brought where the STO describing the thermodynamic limit exhibits a fixed point or limit cycle attractor and shows linear response under perturbations, and examples of STOs with more complicated attractors that do not satisfy linear response.  
  
\subsection{Towards the study of more complicated attractors and behaviors}\label{Sec:TwardsmreComp}
With few exceptions, the situations described in the previous sections can rigorously deal with three cases: small coupling, so that the STO is close to a linear transfer operator; evolution of states close to delta measures for which the STO is close to a finite-dimensional map; or a combination of the two situations. General strategies that deal with genuinely nonlinear infinite-dimensional operators are lacking. 

For example, one wonders if there is a general framework to study STOs with stable fixed point in a regime that is ``far" from linear or a finite-dimensional map, and also if it is  possible to rigorously study bifurcations of STOs where a stable fixed point looses stability. Furthermore, being a nonlinear transformations of an infinite-dimensional space, STOs are expected to have attractors and dynamics more complicated than periodic dynamics. This gives rise to the question  if it is possible to study STOs with multi-dimensional attractors where the dynamics on the attractor and in a neighborhood is amenable to rigorous analysis.

Below we give an example of an innocent looking system of globally coupled maps (perhaps the simplest possible) whose associated STO exhibits very complicated behavior.

\subsubsection{Mean-field coupled rotations} Consider a system of coupled maps where each unit evolves according to a 1D rotation of an angle that depends on the state of all the units via a mean-field. Fix a continuous map $h:\T\rightarrow \R$ and define the system of globally coupled maps
\[
x_i(t+1)=x_i(t)+\frac{1}{N}\sum_{j=1}^Nh(x_i(t)) \mod 1\quad\quad i=1,...N
\]
and the associated  STO
\[
\mc T\mu=\left(f_{\mu}\right)_*\mu
\quad\mbox{where}\quad 
f_\mu(x)=x+\int_\T h(y)d\mu(y)\mod 1.
\] 
This system has a very simple formulation, but as we are going to show, can produce  complicated behavior for the associated  self-consistent transfer operator. 

Consider the one parameter family of rotations $\{R_\theta\}_{\theta\in \T}$ with $R_\theta:\T\rightarrow\T$ and $R_\theta(x)=x+\theta$. Notice that this collection forms a group that acts on the measures in $\mc M(\T)$ as
\[
R_\theta\mu:=(R_\theta)_*\mu.
\] 
Given a measure $\mu$, denote by $C_\mu$ its centralizer, i.e.
\[
C_\mu:=\left\{R_\theta:\, (R_\theta)_*\mu=\mu \right\}.
\] 
If $C_\mu=\{R_0=\Id\}$, we say that $\mu$ has no rotation symmetries. It follows immediately that if $\mu$ has no rotation symmetries, then the orbit of $\mu$ under the $R_\theta$ action, $\{R_\theta\mu\}_{\theta\in \T}$, has a natural 1D torus manifold structure. Let's denote by $\T_\mu\cong \T$ understanding that  $\theta\in\T_\mu$ corresponds to $(R_\theta)_*\mu\in \mc M(\T)$.  

The following proposition is immediate.
\begin{proposition}
Consider $\mu\in \mc M(\T)$ with no rotation symmetries. Then  

i) $\mc T(\T_\mu)\subset \T_\mu$;

ii) $\mc T$ on $\T_\mu$ acts as the map
\begin{equation}\label{Eq:Reduceddynamics}
\mc R_\mu(\theta):=\theta+\int h(y+\theta) d\mu(y)\mod 1.
\end{equation}
\end{proposition}

From the above proposition it follows that the space of measures $\mc M^*$ with no rotation symmetries forms an open subset of $\mc M(\T)$ which, by the above proposition, is foliated into 1D tori on which the dynamics of $\mc T$ has great variability. To see this, take for example an absolutely continuous measure without symmetries having density $\phi$.  Equation \eqref{Eq:Reduceddynamics} becomes
\[
R_\mu(\theta):=\theta+(h*\phi)(\theta)
\]
where $*$ denotes the convolution, and varying $\phi$, i.e. moving from one invariant circle to the next, various maps compatible with the regularity of $h$ can be found\footnote{Whenever $\phi,h\in L^2(\T)$ and admit  a Fourier series expansion,  the image of $\phi\mapsto h*\phi$ can be characterized in terms of the Fourier coefficients of $h$.}. 
%

\subsection{Numerical Studies of STOs} \label{Sec:NumericalStudies}
As the dynamics of STOs can be very complicated to study from a rigorous point of view, computational approaches are crucial to get information on these objects. Here we give two examples where  the attractors of STOs and their bifurcations have been studied numerically. 

\smallskip
In Section \ref{Sec:StabFixedPoints} we reviewed situations where one can prove that  the thermodynamic limit of some coupled systems has a unique equilibrium measure, provided that the coupling strength is sufficiently small. Increasing the coupling strength beyond a certain threshold, one does not expect the picture to persist  and one wonders what kinds of bifurcations the system can undergo. In \cite{selley2021self}  this question has been addressed employing a mix of rigorous arguments and numerical evidence. 

Given $\epsilon\in \R$, consider the parametric family $\{f_\gamma\}_{\gamma\in\R}$ of selfamps of $[0,1]$ 
\[
f_\gamma(x)=\left(2+\epsilon F\left(\frac{1}{\gamma}-2\right)\right)x\mod 1,
\]
and $\hat \gamma:\mc M_1([0,1])\rightarrow \R$  
\[
\hat\gamma(\mu)=\int_{[0,1]} yd\mu(y)
\]
which is the center of mass of $\mu$. For example, if $F(x)=x$, the above becomes
\[
f_{\hat\gamma(\mu)}(x)=\left(2(1-\epsilon)+\frac{\epsilon}{\int yd\mu(y)}\right)x\mod 1
\]
which, for fixed $\mu$, is a $\beta$-transformation giving a perturbation of the doubling map. The self-consistent transfer operator is
\[
\mc T_\epsilon\mu=(f_{\hat\gamma(\mu)})_*\mu.
\]
If $\epsilon=0$, $f_\mu=2x\mod 1$ independently of $\mu$, and the Lebesgue measure is the unique absolutely continuous invariant measure, and  an attracting point for $\mc T_0$. Notice also that  Lebesgue  is a fixed point of $\mc T_\epsilon$ for any value of $\epsilon$.  In \cite{selley2021self}, it has been proven that for $\epsilon>0$ there is another measure absolutely continuous with respect to Lebesgue that is fixed by $\mc T_\epsilon$, and numerical simulations suggest that this measure is a stable fixed point for $\mc T_\epsilon$ while the Lebesgue measure looses stability.

\smallskip
In Section \ref{Sec:LinRespHet} we reviewed numerical evidence showing that in the thermodynamic limit of coupled systems, linear response might fail. This is in disagreement with the chaotic hypothesis of Gallavotti and Cohen claiming that high-dimensional systems are expected to exhibit Axiom A behavior. Motivated by this observation, the work in \cite{wormell2022non} presents an example of an STO fo which numerical evidence suggests the presence of an homoclinic tangency implying robust non-uniformly hyperbolic behaviour.

The system considered there has $\bo x(t)=(x_1(t),...,x_N(t))\in[-1,1]^N$ and
\[
x_i(t+1)=f_{\epsilon \Phi_N(\bo x(t))}(x_i(t+1))\quad\quad i=1,...,N
\]
where $\Phi_N$ is as in \eqref{Eq:MeanfieldpAram}
\[
f_{\alpha}(x) =f(x)+g(\alpha)(1-[f(x)]^2),\quad f(x)=2x-\mbox{sign}(x),\quad g(\alpha)=\frac{3}{16}\cos(8\pi\alpha)
\]
which is a nonlinear perturbation of the doubling map on $[-1,1]$. Notice that with the choice of functions above, the maps $f_\alpha$ are all uniformly expanding with lower bounded uniform expansion. Nonetheless, it is shown  that the STO in the thermodynamic limit has a fixed point that is not attracting, but has some unstable directions that numerical evidence suggests are homoclinic to some stable directions.

\subsection{The thermodynamic limit problem for GCMs}\label{Sec:ThermLimit}
The following question now arises:  to which extent does the thermodynamic limit given by a STO describes the finite dimensional system?  This is a particularly relevant question having in mind applications to  systems composed by a number of units that, although large, has order of magnitude much smaller than e.g. systems from statistical physics. For example, if a macroscopic sample of any gas/solid-state system is composed by $\sim 10^{23}$ molecules, the brain has ``only"  $\sim 10^{10}$ neuronal cells with some substructures (e.g. nuclei and bulbs) counting $\sim 10^3$ neurons. 

In \cite{tanzi2022uniformly}, we provide quantitative estimates for the convergence to the thermodynamic limit in the case where the system of coupled maps is uniformly expanding, and  gives sufficient conditions involving the expansion and interaction strength  ensuring that the limit  approximates  the finite dimensional system for all times up to an error of order $N^{-\gamma}$ with $\gamma<\frac12$, where $N$ is the number of coupled units. As a corollary, one can show that in the limit, the system exhibits \emph{propagation of chaos}\footnote{See Section \ref{Sec:IntPartSyst} for more about propagation of chaos.}.

Furthermore in that paper, globally coupled maps lacking symmetry were introduced. The evolution equations for a system of $N$ \emph{different} coupled maps are
\begin{equation}\label{Eq:DiffCoupMaps}
x_i(t+1)= F_i(x_i(t);\, \hat {\bo x}_i(t)):=f_i(x_i(t))+\frac{1}{N}\sum_{j=1}^Nh_{ij}(x_i(t),x_j(t)).
\end{equation}
where $\hat {\bo x}_i(t)$ gathers all coordinates but the $i$-th one. To describe the behavior of this system when $N$ is finite but large, the STO acting on $\mc M(\T^N)$ defined in \eqref{Eq:SelfConsistNoSymm} was introduced. The heuristic behind this definition is that, considering a product probability measure $\mu=\mu_1\otimes...\otimes\mu_N$, concentration results suggest that when  $N$ is large and $(x_1,...,x_N)$ is sampled according to $\mu$, the average in \eqref{Eq:DiffCoupMaps} can be approximated by  its expectation with respect to $\mu$ with high probability
\begin{equation}\label{Eq:Concentrationofmeasure}
\frac{1}{N}\sum_{j=1}^Nh_{ij}(x_i,x_j)\approx \frac{1}{N}\sum_{j=1}^N\int_\T h_{ij}(x_i,\,y)d\mu_j(y).
\end{equation}

To prove the result above, it has been showed that these globally coupled maps preserve a class of measures that are close enough to being  products, so that they satisfy the usual concentration properties of independent bounded random variables, and with respect to which the  approximation in \eqref{Eq:Concentrationofmeasure}  holds with high probability, thus allowing to greatly simplify the evolution equation.  Roughly speaking, this implies that picking an initial condition with respect to a measure from this class,  with high probability, the evolution is undistinguishable from application of the STO up to a small error. Crucial to this argument is the study of the evolution of conditional measures on non-invariant foliations with leaves along the coordinate directions.

\section{Mean-field models in continuous time and other topics}\label{Sec:Othertopics}
The objective of this section is to give some pointers to other branches of the literature on the study of mean-field interacting systems beyond the study of coupled maps. Some of the topics we mention are established fields of research, and our exposition is going to be very superficial with no pretense of completeness.  
\subsection{Coupled systems in continuous time}

\subsubsection{Mean-field interacting particle systems and propagation of chaos}\label{Sec:IntPartSyst} As a prototypical example of mean-field models in continuous time consider the Markov process $(X^N_1(t),...,X^N_N(t))\in \R^{Nd}$ describing $N$ identical entities coupled through a mean-field  plus noise:
\begin{equation}\label{Eq:IntSyst}
dX^N_i(t)=\frac{1}{N}\sum_{j=1}^N h(X_i^N(t),\,X_j^N(t))\,dt  + dW_i(t)\quad\quad i=1,...,N
\end{equation}
where $h:\R^d\times\R^d\rightarrow\R^d$ is a coupling function and $(W_i(t))_i$ are $N$ independent Brownian motions. It is well known  that if $h$ is Lipschitz and one pick an initial condition such that $(X_i^N(0))_{i=1}^N$ are i.i.d., then the SDE admits a solution. Other important models of interacting particles are deterministic (there is no $dW_i(t)$ in the equations) or $h$ is singular\footnote{For example, when $h$ comes from Coulomb interactions. Studying equations \eqref{Eq:IntSyst} in this case has been a long standing problem that was recently solved in \cite{serfaty2020mean}.}. 

The main objective is to study the above system when $N\rightarrow \infty$. To this end, one can define
\[
\mu^N_t=\frac{1}{N}\sum_{j=1}^N\delta_{X_i^N(t)}
\]
and \eqref{Eq:IntSyst} becomes
\[
dX_i^N(t)=H(X_i^N(t),\mu^N_t)dt+dW_i(t)\quad\quad i=1,...,N.
\]
with
\[
 H(X_i^N(t),\mu^N(t)):=\int h(X_i^N(t),y)\,d\mu^N_t(y).
\]
This is the continuous time analogue of the discrete time evolution given in \eqref{Eq:Selfcons2} in Sect. \ref{Sec:ViolLN}.

The weak coupling among the components and the choice of initial condition with i.i.d. components suggest that in the limit  $N\rightarrow \infty$, for every $k\in \N$ and $t\ge 0$,
\begin{equation}\label{Eq:PropofChaos}
(X_1^N(t),...,X_k^N(t))\rightarrow (\bar X_1(t),...,\bar X_k(t))
\end{equation}
in distribution, where $(\bar X_i(t))_{i=1}^N$ are i.i.d processes solving 
\begin{equation}\label{Eq:Limit}
d\bar X(t)=H(\bar X(t),\bar \mu_t)dt+dW(t)
\end{equation}
where $\bar \mu_t$ is the distribution of $\bar X(t)$. These equations are also referred to as \emph{McKean-Vlasov equations}. When \eqref{Eq:PropofChaos} holds, the system is said to exhibit \emph{propagation of chaos}. Notice that the self-consistent transfer operator is the  discrete time analogue of the generator of \eqref{Eq:Limit}. The presence of noise and the exchangeability of the system (full permutation symmetry) are some fundamental ingredients to prove propagation of chaos in the setup above. 
 
In general, to make the above picture rigorous one has to check that
\begin{itemize}
\item[i.] the  system of coupled SDEs in \eqref{Eq:IntSyst} is well posed\footnote{If $h$ has singularities this can be nontrivial.};
\item[ii.] the self-consistent SDE in \eqref{Eq:Limit} describing the candidate limit is well-posed;
\item[iii.] the limit in \eqref{Eq:PropofChaos} holds.
\end{itemize} 

Various approaches to prove \eqref{Eq:PropofChaos} exist among which coupling methods were  the first to be employed, while large deviation estimates and entropy bounds are among the most recent. The full permutation symmetry of the system is central to the above analyses. 

Notice that in discrete time: point i. corresponds to the definition of the finite dimensional map, and ii. corresponds to existence of the mean-field map  so there is nothing to prove; an effort has to be put in iii., especially to gain explicit rates of convergence (and this is the topic of Section \ref{Sec:ThermLimit}); while most of the analysis of STOs presented in previous sections corresponds in the continuous time setup to the study of the solutions of \eqref{Eq:Limit}.

For more on interacting particle systems: \cite{sznitman1991topics}  is a classical reference on propagation of chaos;    \cite{spohn2012large} and \cite{cercignani2013mathematical}  give accounts on the study of interacting systems in statistical mechanics; \cite{giacomin2019stochastic} contains a nice introduction to the topic together with a discussion of applications and some recent contributions;  \cite{chaintron2021propagation} is a very thorough review containing an account of various setups of interacting particle systems, detailed definition(s) of propagation of chaos and different approaches to prove it.

\subsubsection{Coupled Oscillators}\label{Sec:CoupledOscill}
Oscillating systems are ubiquitous in nature and artificial systems (see Section \ref{Sec:ApplicationofMeanField}) and often interact with one another. An example of coupled oscillators is  given by the system of $N$ coupled differential equations
\[
\frac{d\theta_i}{dt}=\omega_i+\frac{K}N\sum_{j=1}^Nh_{ij}(\theta_j-\theta_i)\quad\quad i=1,...,N
\]
where  $\omega_i$ are called \emph{natural frequencies} and give the angular velocities at which the oscillators would rotate if uncoupled. They are in general all different and are usually assumed to be real i.i.d. random variables having distribution with some density $g$.

This  model originated as a simplification of the dynamics of weakly coupled (almost identical ODEs) having an attracting limit cycle (see \cite{winfree1967biological}). It was first intensely  studied by Kuramoto in the case where $h_{ij}$ is the sine function yielding what is known as the Kuramoto model \cite{kuramoto1984chemical}:
\begin{equation}\label{Eq:Kuramoto}
\frac{d\theta_i}{dt}=\omega_i+\frac{K}N\sum_{j=1}^Nh_{ij}(\theta_j-\theta_i) \quad\quad i=1,...,N.
\end{equation}

The main starting observation is that when the coupling strength $K\ge 0$ is small,  the difference in the natural frequencies causes the oscillators to spread and the system appears to be disordered. When $K$ is increased above a certain threshold the system synchronizes. 

A first attempt to study this phenomenon involved the introduction of variables $r$ and $\psi$ defined as
\[
r(t)e^{i\psi(t)}=\frac{1}{N}\sum_{j=1}^Ne^{i\theta_j(t)}
\] with respect to which the equations in \eqref{Eq:Kuramoto} can be written (after a change of variables) as
\[
\frac{d\theta_i}{dt}=\omega_i-Kr\theta_i\quad\quad i=1,...,N.
\]
Notice that $r=0$ implies a ``disordered" distribution of the angles $\theta_i$, while $r=1$ implies $\theta_1=\theta_2=....=\theta_N$, i.e. a fully synchronized state, therefore 
this parameter gives important information on the state of the system.

The evolution of $r(t)$ has been rigorously studied in \cite{strogatz1991stability} for $N\rightarrow \infty$. In this limit, the state of the system at time $t$ was assumed to be described by the collection of densities $\rho_{\omega,t}(\theta):=\rho(\theta,\omega,t)$ -- this corresponds to the state of the GCM in the thermodynamic limit --, i.e. the density of the distribution of the oscillators having natural frequency $\omega$ at time $t$, with the function $\rho$ satisfying a PDE originating from a continuity equation -- which gives an evolution law  corresponding to the STO in the GCM setup.

The study of the Kuramoto model generated a large body of works. We direct the reader to the excellent paper  \cite{strogatz2000kuramoto}. For other reviews see also: \cite{acebron2005kuramoto}; \cite{pikovsky2015dynamics} for   generalizations  and other models of coupled oscillators; \cite{rodrigues2016kuramoto} for  works on oscillators coupled on various types of networks.


\subsubsection{Chimeras}\label{subsec:chimeras} Roughly speaking, a chimera is a state of a system of coupled units  where part of the units exhibit coherent motion, e.g. are synchronized, while another part exhibit erratic incoherent behavior (this is analogous to the partially ordered phase described in Section \ref{Sec:GlobCouplVarBehav}). Most surprisingly, chimeras were observed early on in systems of coupled oscillators having full permutation symmetries suggesting the presence of symmetry breaking \cite{abrams2004chimera,abrams2006chimera}.

There is no general consensus of what constitutes or not a chimera. A review of different characterizations can be found in \cite{haugland2021changing}, while \cite{bick2016chaotic} proposed a rigorous mathematical definition. In some cases chimera states are expected to arise as stationary states for the system \cite{martens2016basins} while in others arise only as long transients \cite{wolfrum2011chimera}.

\cite{panaggio2015chimera,yao2016chimera,bera2017chimera} are reviews focusing on the phenomenology of chimeras, while \cite{omel2018mathematics} focuses on their mathematical study.

\subsubsection{Coupled Systems and their Symmetries} 

Given $X\subset\R^N$ and a vector field $f:X\rightarrow X$, the ODE
\[
\frac{dx}{dt}=f(x)
\] 
has the linear transformation $\gamma:\R^N\rightarrow \R^N$ as a symmetry if $f(\gamma x)=\gamma x$. Permutation symmetries are an example where one or more coordinates can be swapped without changing the vector field\footnote{\eqref{Eq:IntSyst} without the noise term, gives an ODE with full permutation symmetry.}.

When the vector field describes the continuous time evolution of units coupled on a graph, it is likely that the symmetries of the graph correspond to symmetries of the dynamics, and therefore contain information on the time evolution \cite{golubitsky2003symmetry,golubitsky2003symmetry,field2009symmetry}. In fact, they can explain synchronization and coherence patterns as well as the bifurcations that these patterns undergo when the parameters of the system change. This is a notable example of how the interacting structure can influence the dynamic. 

 For a review see \cite{golubitsky2015recent}. See also \cite{aguiar2011dynamics,rink2014coupled,rink2015coupled} and references therein. For a study of the role of symmetry in a system of identical coupled oscillators see \cite{ashwin1992dynamics}. For a study of the role of symmetries on mean-field limits see \cite{bick2021mean}.

\subsection{Different types of couplings}

Interactions shape the dynamics of complex systems and can produce behavior drastically different from that of the local dynamics.  It was noted in the previous section that the coupling structure (who is connected to whom) influences the resulting dynamic. The particular form of the coupling (e.g. the function $h$ in \eqref{Eq:EqsCoupledSyst}) has also an important role (see e.g.  \cite{stankovski2017coupling}).

In contrast with what we have presented so far, in this section we review works where the coupling changes with time (\emph{adaptive networks}); and when the coupling can arise among multiple units (\emph{higher-order networks}), rather than from pairwise interactions.

\subsubsection{Adaptive networks} In this type of networks, the interaction among the units changes according to the internal state of the system. For example, some links in the graph of interactions could be severed or added with time, or more generally the coupling strength among different nodes could be increased or decreased. These networks capture many phenomena in real-world systems, a prime example is plasticity between neurons by which signal transmission at synapses is strengthened or weakened to modify the dynamic (plasticity is at the base of development, learning, and memory). 

For a mean-field reduction approach to adaptive networks related to plasticity see \cite{duchet2022mean} and references therein. For an application to networks on power grids see \cite{berner2021adaptive}.
For many more  on adaptive networks and examples arising in real-world systems see \cite{gross2008adaptive,gross2009adaptive,sayed2014adaptive}, for an earlier reference see \cite{sutton1981toward}.

\subsubsection{Higher-order networks}

All the  coupled systems considered so far have been characterized by pairwise additive interactions, meaning that the interaction term is given by the sum of all the interactions between a node and each other node in the network. In contrast, systems coupled in higher-order networks (\cite{bianconi2021higher,bick2021higher}) have interactions terms where the interaction can also be among 3 or more units at the same time. For example, one can imagine a situation where the interaction strength between two nodes is modulated by a third node in the network. In this case,  the interaction term depends on the coordinates of all three nodes. 

As an example, the Kuramoto model discussed in Section \ref{Sec:CoupledOscill} arises as a first order approximation of a system of interacting limit cycles with only  pairwise interactions; in contrast,  \cite{bick2016chaos} contains a derivation of the higher orders where the interaction terms depend on multiple oscillators.  Higher-order networks have recently shown to arise also as a result of the choice of coordinates \cite{nijholt2022emergent}.  

Mean-field coupled higher-order networks and their thermodynamic limits have been studied in \cite{bick2022multi,bick2022phase,gkogkas2022graphop}.

 \subsection{Mean-field coupled models as models of real-world systems}\label{Sec:ApplicationofMeanField}
In this section we give some indications to reviews and some selected works on modeling real-world systems via mean-fields.

\smallskip
Globally coupled oscillators and maps arise as models of several physical systems. Some lists of applications can be found in the introductions to \cite{kaneko1989pattern,kaneko1991globally,strogatz2000kuramoto}. Among these, one finds Josephson junctions array, charge density waves, nonlinear optics, coupled lasers, and microwave oscillators. 

\smallskip
Mean-field coupled maps and flows have received particular attention for their ability to reproduce behavior observed in systems of biological origin. This was very early on noted, among others, by Kaneko   \cite{kaneko1994relevance} (see also the more recent \cite{kaneko2015globally}) who reviewed several biological systems in which coherent structures like those described in Section \ref{Sec:GlobCouplVarBehav} were observed to arise as the result of the interaction of many components. Another feature of GCM that recalls the functioning of some biological systems is the presence of a great variety of attractors that the system can visit under perturbation due to external forcing. This characteristics can allow for some computational mechanisms: as the external factors change, orbits are sent to a different attractor that encodes a particular stimulus or some features of the stimulus. 

\smallskip
The use of mean-field models to simulate the behavior of globally coupled neurons has a long list of contributions. For what concerns map-based models of neurons (i.e. models in discrete time) a review is given in  \cite{ibarz2011map} which contains a list of studies on globally coupled maps describing the evolutions of ensembles of neurons.  Here we mention \emph{Rulkov Maps} \cite{rulkov2002modeling} which  model  chaotic bursting, a firing  pattern where stretches of high-frequency spiking are alternated (in an erratic fashion) with stretches where the neuron is at rest.  A system of coupled Rulkov Maps is given by equations (for $i=1,...,N$)
\begin{equation*}
\left\{\begin{array}{lcl}
x_i(t+1)&=&\frac{\alpha}{1+x_i(t)^2}+y_i(t)+\frac{\epsilon}{N}\sum_{j=1}^Nx_j(t)\\
y_i(t+1)&=&y_i(t)-\sigma x_i(t)-\beta
\end{array}\right.
 \end{equation*}
with parameters $\alpha\approx 4.2$ and $\sigma=\beta\approx0.001$. 
 The fast variable $x_i$ describes the membrane voltage of the neuron, while the slow variable $y_i$ describes an internal variable that is responsible for the switching between resting and chaotic phase. Mean-field coupled Rulkov Maps and their synchronization are studied in\cite{rulkov2001regularization}. 

\smallskip
A review of models  of coupled oscillators that capture some aspects of  neuron dynamics can be found in \cite{bick2020understanding}. Here we mention Ermentrout and Koppell's \emph{Theta Model} which is a continuous time 1D model of tonic spiking activity in neurons. A single neuron is described by the ODE on the unit circle
\[
\frac{d}{dt}\theta=1-\cos\theta+(1+\cos\theta)(r+I(t))
\]
where $r<0$ is a resting potential, and $I(t)$ is the current coming to the neuron at time $t$.
When $I(t)<|r|$, the system has an attracting fixed point, and the dynamics is at ``rest". At $I(t)=|r|$ the system undergoes a saddle node bifurcation, and for $I(t)>|r|$ orbits rotate around the circle corresponding to a ``tonic spiking" phase. Coupling Theta Models one obtains the equations  
\[
\frac{d}{dt}\theta_i(t)=1-\cos\theta_i+(1+\cos\theta_i)(r+I_i(t))\quad\quad i=1,...,N
\]
where $I_i(t)$ now depends on $\{\theta_j(t)\}_{j=1}^N$, and for example (for instantaneous synapses) are equal to
\[
I_i(t)=\frac1N\sum_{j=1}^N(1-\cos\theta_j(t))^2.
\] 
A study of the dynamics is undertaken using approaches similar to those described in Section \ref{Sec:CoupledOscill} for the Kuramoto model  \cite{kotani2014population,laing2018dynamics}. 

\smallskip
The Kuramoto model  can also arise as a phase reduction model of coupled neurons \cite{izhikevich1999weakly}. Other  works study mean-field coupled models of integrate and fire neurons (called \emph{population density models} in the computational neuroscience literature) and can reproduce some characteristics  oscillations recorded in brain activity that are known as \emph{rythms} \cite{nykamp2000population,haskell2001population}. For a review of population density models see \cite{Brunel2015}. 

\smallskip
Mean field coupled models have been proposed also to reproduce the behavior of gene networks. A  general strategy to obtain mean-field models for gene regulatory networks has been proposed in \cite{andrecut2006mean}.  An interesting mean-field model we mention is the one  where the expression of a group of genes is regulated by a common repressor field. This generated a simplified model of mean-field coupled \emph{degrade and fire oscillators}  with interesting features, like clustering, and which is amenable to rigorous analysis and classification of the periodic attractors and their basins \cite{fernandez2014typical,fernandez2018revisiting}. 

\smallskip
Related to applications is also the problem of recovering models for coupled systems from observational data. This is a particularly hard task for mean-field coupled systems where the very small size of the interactions hinders reconstruction of the connections among units via model based methods, and the erratic dynamics makes  model free estimation (e.g. correlation analysis) ineffective. For some  contributions to this problem that specifically address mean-field coupled systems and the issues that arise in this set-up see \cite{eroglu2020revealing}.

\smallskip 
Mean-field models have also been extensively studied in game theory see e.g. \cite{bender2003self,lasry2007mean,carmona2013probabilistic}; in opinion models \cite{kacperski1999opinion}; and social sciences \cite{contucci2008phase,gallo2009parameter}.

\appendix

\section{Some proofs from Section \ref{Sec:StabFixedPoints}}\label{App:Comput}
\begin{proof}[Proof of Lemma \ref{Lem:InvarianceBallLip}]
For a fixed $\phi\in B_L$, $g_\phi(x)=x+\epsilon\int_\T h(x,y)\phi(y)dy$. For $|\epsilon|$ small enough, $g_\phi$ is a diffeomorphism and 
\begin{align*}
\mc L_{\epsilon,\phi}(\phi)(x)= \frac{\phi}{|g_\phi'|}\circ g_\phi^{-1}(x)
\end{align*}
A computation shows that above is Lipschitz, with Lipschitz constant $(1+O(\epsilon))(L+O(\epsilon))$ with $O(\epsilon)$ uniform in $\phi$. Another computation shows that $\mc P$, defined in \eqref{Eq:TransfOpDoubMap}, halves Lipschitz constants. Therefore $\mc P\mc L_{\epsilon,\phi}(\phi)$ has Lipschitz constant 
\[
\frac{1}{2}(1+O(\epsilon))(L+O(\epsilon))
\]
and the result follows. 
\end{proof}

\begin{proof}[Proof of Lemma \ref{Lem:ContC0Top}]
$\mc P$ is evidently continuous. For $\phi\mapsto \mc L_{\epsilon,\phi}\phi$, triangle inequality implies
\begin{align}
|\mc L_{\epsilon,\phi_1}\phi_1(x)-\mc L_{\epsilon,\phi_2}\phi_2(x)| &=  \left|\frac{\phi_1}{|g_{\phi_1}'|}\circ g_{\phi_1}^{-1}(x)-\frac{\phi_2}{|g_{\phi_2}'|}\circ g_{\phi_2}^{-1}(x) \right|\nonumber \\
&\le \left|\frac{\phi_1}{|g_{\phi_1}'|}\circ g_{\phi_1}^{-1}(x)-\frac{\phi_2}{|g_{\phi_1}'|}\circ g_{\phi_1}^{-1}(x)\right|+\label{Eq:EstTriang1}\\
&\quad\quad+ \left|\frac{\phi_2}{|g_{\phi_1}'|}\circ g_{\phi_1}^{-1}(x)-\frac{\phi_2}{|g_{\phi_1}'|}\circ g_{\phi_2}^{-1}(x)\right|\label{Eq:EstTriang2}\\
&\quad\quad+ \left|\frac{\phi_2}{|g_{\phi_1}'|}\circ g_{\phi_2}^{-1}(x)-\frac{\phi_2}{|g_{\phi_2}'|}\circ g_{\phi_2}^{-1}(x)\right|.\label{Eq:EstTriang3}
\end{align}
For the  term in \eqref{Eq:EstTriang1}
\[
\left|\frac{\phi_1}{|g_{\phi_1}'|}\circ g_{\phi_1}^{-1}(x)-\frac{\phi_2}{|g_{\phi_2}'|}\circ g_{\phi_2}^{-1}(x) \right|\le(1+O(\epsilon))\|\phi_1(y)-\phi_2(y)\|_{C^0}
\]
where $O(\epsilon)$ is uniform in $\phi$.
To bound the  term in \eqref{Eq:EstTriang2}, notice that 
\begin{align*}
|g_{\phi_1}^{-1}(x)-g_{\phi_2}^{-1}(x)| &\le d_{C^0}(g_{\phi_1},g_{\phi_2})  \sup_{y\in \T}|g_{\phi_1}'(y)|\\
&\le O(\epsilon)\|\phi_1-\phi_2\|_{C^0}
\end{align*}
where $O(\epsilon)$ is a constant that depends only on the derivatives of $h$ and $\epsilon$, and is uniformly bounded when $\epsilon$ is bounded. This implies that
\begin{equation}\label{Eq:BoundC0norm}
\left|\frac{\phi_2}{|g_{\phi_1}'|}\circ g_{\phi_1}^{-1}(x)-\frac{\phi_2}{|g_{\phi_1}'|}\circ g_{\phi_2}^{-1}(x)\right|\le \left| \frac{\phi_2}{|g_{\phi_1}'|}\right|_{Lip} O(\epsilon)\, \|\phi_1-\phi_2\|_{C^0}.
\end{equation}
For the term in \eqref{Eq:EstTriang1}, 
\begin{equation}\label{Eq:BOundjacobian}
\left|\frac{1}{|g_{\phi_1}'|}(x)-\frac{1}{|g_{\phi_2}'|}(x)\right|\le K_2\left||g_{\phi_1}'|-|g_{\phi_2}'|\right| \le K_3 \|\phi_1-\phi_2\|_{C^0}
\end{equation}
with $K_2$ and $K_3$ independent of $x$, $\phi_1$, and $\phi_2$; in the first inequality we used that $|g_{\phi_i}'|(x)$ is uniformly bounded away from zero in both $x$ and $\phi_i$ while in the second we used 
\[
|g_{\phi_1}'|-|g_{\phi_2}'|\le 2 \epsilon \int |\partial_1 h| |\phi_1-\phi_2|\le O(\epsilon) \|\phi_1-\phi_2\|_{C^0}.
\]

Putting together all the inequalities above the lemma is proved.
\end{proof}

\begin{proof}[Proof of Lemma \ref{Lem:Continuityrelation}]
Since $\mc P$ is bounded in the $\|\cdot\|_{C^1}$ norm, it is enough to prove that $\exists K'\ge 0$ such that for all $\phi\in\mc B_{L_1,L_2}$
\[
\|( \mc L_{\epsilon,\phi} -\mc L_{\epsilon,\phi_*}) \phi_*\|_{C_1}\le K' \epsilon\cdot \|\phi-\phi_*\|_{C_1}.
\]
Now, $\|( \mc L_{\epsilon,\phi} -\mc L_{\epsilon,\phi_*}) \phi_*\|_{C_0}$ has been already bounded in \eqref{Eq:BoundC0norm}. We proceed with a bound of the first derivative
\begin{align}
\left|[( \mc L_{\epsilon,\phi} -\mc L_{\epsilon,\phi_*}) \phi_*]'\right| &= \left| \left[ \frac{\phi_*}{|g_{\phi}'|}\circ g_{\phi}^{-1}(x)-\frac{\phi_*}{|g_{\phi_*}'|}\circ g_{\phi_*}^{-1}(x)\right]'\right|\\
&\le \left|\frac{\phi_*'}{|g_{\phi}'|^2}\circ g_{\phi}^{-1}(x)-\frac{\phi_*'}{|g_{\phi_*}'|^2}\circ g_{\phi_*}^{-1}(x)\right| +\label{Eq:EstLine1}\\
&\quad\quad + \left|-\frac{\phi_*g_{\phi}''}{|g_{\phi}'|^3}\circ g_{\phi}^{-1}(x)+\frac{\phi_*g_{\phi_*}''}{|g_{\phi_*}'|^3}\circ g_{\phi_*}^{-1}(x)\right| \label{Eq:EstLine2}
\end{align}
For \eqref{Eq:EstLine1}
\begin{align*}
\left|\frac{\phi_*'}{|g_{\phi}'|^2}\circ g_{\phi}^{-1}(x)-\frac{\phi_*'}{|g_{\phi_*}'|^2}\circ g_{\phi_*}^{-1}(x)\right| &\le \left|\frac{\phi_*'}{|g_{\phi}'|^2}\circ g_{\phi}^{-1}(x)-\frac{\phi_*'}{|g_{\phi}'|^2}\circ g_{\phi_*}^{-1}(x)\right| +\\
&\quad\quad+\left|\frac{\phi_*'}{|g_{\phi}'|^2}\circ g_{\phi_*}^{-1}(x)-\frac{\phi_*'}{|g_{\phi_*}'|^2}\circ g_{\phi_*}^{-1}(x)\right| \\
&\le \left|\frac{\phi_*'}{|g_{\phi}'|^2}\right|_{Lip} O(\epsilon)\|\phi-\phi_*\|_{C^0}+\\
&\quad\quad+\|\phi_*\|_{C^1}O(\epsilon)\|\phi-\phi_*\|_{C^0}\\
&\le O(\epsilon)\|\phi-\phi_*\|_{C^0}
\end{align*}
where the bound on the first term follows as \eqref{Eq:BoundC0norm}, and for the second term
\[
\left|\frac{\phi_*'}{|g_{\phi}'|^2}\circ g_{\phi_*}^{-1}(x)-\frac{\phi_*'}{|g_{\phi_*}'|^2}\circ g_{\phi_*}^{-1}(x)\right|\le \|\phi_*\|_{C^1}\,O(\epsilon)\,\||g_{\phi_*}'|^{-2}-|g_{\phi}'|^{-2}\|_{C^0}
\]
and the above can be bound with computations analogous to \eqref{Eq:BOundjacobian}.

Similar estimates imply that the expression in \eqref{Eq:EstLine2} can be bounded by $O(\epsilon)\|\phi-\phi_*\|_{C^0} $.
\end{proof}
\begin{proof}[Proof of Lemma \ref{Lem:COnes}]
\begin{align*}
\frac{\phi}{|g_\phi'|}\circ g_\phi^{-1}(x)\left(\frac{\phi}{|g_\phi'|}\circ g_\phi^{-1}(y)\right)^{-1}&=\frac{\phi\circ g_\phi^{-1}(x)}{\phi\circ g_\phi^{-1}(y)}\cdot \frac{|g_\phi'|\circ g_\phi^{-1}(y)}{|g_\phi'|\circ g_\phi^{-1}(x)}
\end{align*}
and the key estimates are
\[
\frac{\phi\circ g_\phi^{-1}(x)}{\phi\circ g_\phi^{-1}(y)}\le e^{a| g_\phi^{-1}(x)- g_\phi^{-1}(y)|}\le e^{a(1+ O(\epsilon))|x-y|}
\]
and
\begin{align*}
\frac{|g_\phi'|\circ g_\phi^{-1}(y)}{|g_\phi'|\circ g_\phi^{-1}(x)}&=\frac{1+\epsilon\int \partial_1h(g_\phi^{-1}(x),z)\phi(z)dz}{1+\epsilon\int \partial_1h(g_\phi^{-1}(y),z)\phi(z)dz}\\
&\le 1+\epsilon\frac{\int |\partial_1^2h|_\infty \phi(z)dz}{1+\epsilon\int \partial_1h(g_\phi^{-1}(y),z)\phi(z)dz}\cdot  | g_\phi^{-1}(x)- g_\phi^{-1}(y)|\\
&\le 1+ O(\epsilon)\cdot (1+\mc O(\epsilon))|x-y|\\
&\le e^{ O(\epsilon)|x-y|}.
\end{align*}
\end{proof}

\bibliographystyle{amsalpha}
\bibliography{Bibliography}

\end{document}